# Enumerative geometry of $K3$ surfaces
## by Thomas Dedieu

## 1 – Introduction

The aim of these notes is to explain various enumerative results about $K3$ surfaces without assuming familiarity with Gromov–Witten theory; in fact, they represent an attempt on my part to understand the meaning of these results in classical terms. The enumerative results in question are due to Beauville, Bryan and Leung, Pandharipande, Maulik, Thomas, and others, and confirm conjectures made by Yau–Zaslow, Göttsche, and Katz–Klemm–Vafa. They are listed in Paragraph (1.1) below.

They fall in three categories: (i) results that don't really need Gromov–Witten theory at all, either to be formulated or to be proved; (ii) results that may be formulated without Gromov–Witten theory, but which necessitate techniques from this theory to be proven (as far as current knowledge goes); (iii) results which require an understanding of Gromov–Witten theory to be fully appreciated. The text provides the necessary background in Gromov–Witten theory in order to understand results from all these three categories. In particular we shall give a brief idea of what Gromov–Witten invariants are in Paragraph (1.2).

**(1.1) Contents description.** In Section 2 we state a formula of Yau and Zaslow giving the number of rational curves in a primitive linear system $|L|$ on a $K3$ surface $S$, and give its proof by Beauville following the strategy suggested by Yau–Zaslow, using the universal compactified Jacobian; this falls in the above mentioned category (i). In this formula, each rational curve is counted with multiplicity equal to the topological Euler number of its compactified Jacobian, which depends only on the singularities of the curve. We also present results due to Fantechi, Göttsche, and van Straten which show that this is also the multiplicity of the origin in a semi-universal deformation space of the singularities of the curve.

The Yau–Zaslow formula is generalized in Section 3 to a formula giving the number of genus $g$ curves in the system $|L|$ passing through $g$ general points on $S$, for arbitrary $g$; this generalized formula had been conjectured by Göttsche. We give an outline of its proof by Bryan and Leung, by degeneration to an elliptic $K3$ surface, with as much details as the scope of these notes and the ability of the author permit. This proof requires the formulation of the result in terms of some special version of Gromov–Witten invariants specifically designed for algebraic $K3$ surfaces (see Section 3.1), and relies among other things on a multiple cover formula for nodal rational curves. This result thus falls in category (ii).

Essentially all remaining results fall in category (iii). The goal of Section 4 is to explain the extension of the Yau–Zaslow formula to a formula for non-primitive linear systems, which has been proven by Klemm–Maulik–Pandharipande–Scheidegger (this proof is streamlined in Section 6.3). This features the Aspinwall–Morrison multiple cover formula, and its application to define corrected Gromov–Witten invariants known as BPS states numbers. I also discuss other degenerate contributions to these invariants besides plain multiple curves, striving to sort out the relation between the number given by the formula and the actual number of integral rational curves.

In Section 5 we introduce various generalizations of the Gromov–Witten invariants on $K3$ surfaces considered so far, in particular the Hodge integrals and their BPS states version. The main theme in this section is the relation between between Gromov–Witten integrals on $K3$ surfaces and curve counting Gromov–Witten invariants on threefolds. For instance, we discuss the Katz–Klemm–Vafa formula, proved by Pandharipande and Thomas, which may be seen as computing the excess contribution of an algebraic $K3$ surface to the Gromov–Witten invariants of any fibered threefold in which it appears as a fibre. Some of the material presented in this section is necessary for the proof (discussed in Section 6) of the Yau–Zaslow formula for non primitive classes (stated in Section 4).

In Section 6 we introduce Noether–Lefschetz numbers for families of lattice-polarized $K3$ surfaces, and state a theorem due to Maulik and Pandharipande which shows, on a threefold fibered in lattice-polarized $K3$ surfaces, how these Noether–Lefschetz numbers give an explicit relation between Gromov–Witten invariants of the threefold and of the $K3$ fibres. Finally, we discuss the application of this formula to the proof of the Yau–Zaslow formula for non-primitive linear systems. It involves a mirror symmetry theorem that enables the computation of Gromov–Witten invariants of anticanonical sections of toric 4-manifolds, as well as modularity results for Noether–Lefschetz numbers following from the work of Borcherds and Kudla–Millson; the latter enable the computation of all Noether–Lefschetz numbers of the family of lattice-polarized $K3$ surfaces considered in the proof, viz. the so-called $STU$ model.

**(1.2) Gromov–Witten theory.** Let $X$ be a compact complex manifold. In Gromov–Witten theory, one views genus $g$ curves on $X$ as *stable maps*, i.e., morphisms $f : C \to X$ where $C$ is a connected, nodal, curve of arithmetic genus $g$, such that there are only finitely many automorphisms $\phi$ of $C$ satisfying the identity $f \circ \phi = f$. The latter condition is called the *stability condition*, and amounts to the requirement that each irreducible component of arithmetic genus 0 (resp., 1) of $C$ which is contracted by $f$ carries at least 3 (resp., 1) *special points*, i.e., either intersection points with other irreducible components of $C$ or, if relevant, marked points. An integral embedded curve $C \subseteq X$ is then encoded as the map $f : \bar{C} \to X$ obtained by composing the normalization of $C$ with its embedding in $X$.

Stable maps have moduli spaces that are compact. Let $g$ and $k$ be non-negative integers, and $\beta$ be a homology class in $H_2(X, \mathbf{Z})$. We denote by $\overline{M}_{g,k}(X, \beta)$ the moduli space of genus $g$ stable maps $f : C \to X$, with $k$ marked points on $C$, such that $[f_*(C)] = \beta$. It is compact, has expected dimension

$$(1.2.1) \qquad \operatorname{expdim}\bigl(\overline{M}_{g,k}(X,\beta)\bigr) = (\dim(X) - 3)(1 - g) - K_X \cdot \beta + k,{}^1$$

and all its irreducible components have dimension larger than or equal to the expected dimension. The moduli space $\overline{M}_{g,k}(X, \beta)$ has in general many irreducible components, including some parametrizing stable maps that one would arguably rather discard if one is willing to enumerate genuine curves on $X$: for instance, note that the contraction of a curve $C$ of genus $g \geqslant 2$ to a point of $X$ is a stable map. It turns out nevertheless that considering stable maps instead of embedded curves provides much better behaved (e.g., deformation invariant) enumerative numbers. It is possible to consider different devices instead of stable maps in order to obtain well-behaved enumerative numbers, each such possibility coming with its own specific pros and cons; about this, see the enlightening survey [32].

Now, Gromov–Witten invariants are[2] intersection numbers

$$(1.2.2) \qquad \int_{[\overline{M}_{g,k}(X,\beta)]^{\mathrm{vir}}} \operatorname{ev}_1^*(\gamma_1) \cup \ldots \cup \operatorname{ev}_k^*(\gamma_k),$$

where: $\gamma_1, \ldots, \gamma_k$ are cohomology classes in $H^*(X, \mathbf{Z})$; for all $i = 1, \ldots, k$, $\operatorname{ev}_i$ is the *evaluation map* at the $i$-th marked point, i.e., the map $\overline{M}_{g,k}(X, \beta) \to X$ defined by $[f : C \to X, x_1, \ldots, x_k] \mapsto f(x_i)$; $[\overline{M}_{g,k}(X, \beta)]^{\mathrm{vir}}$ is a rational homology class in $H_{2\,\mathrm{vdim}}\bigl(\overline{M}_{g,k}(X, \beta), \mathbf{Q}\bigr)$, called the *virtual fundamental class*, where vdim is the expected (a.k.a. virtual) dimension of $\overline{M}_{g,k}(X, \beta)$. The virtual fundamental class is the usual fundamental class when the moduli space $\overline{M}_{g,k}(X, \beta)$ has the expected dimension, but in general it is given by an excess formula (it is the top Chern class of the obstruction bundle when $\overline{M}_{g,k}(X, \beta)$ is non-singular). Thus, the virtual fundamental class allows the computation of the intersection product $\operatorname{ev}_1^*(\gamma_1) \cup \ldots \cup \operatorname{ev}_k^*(\gamma_k)$ on $\overline{M}_{g,k}(X, \beta)$ "as if the latter were equidimensional of the expected dimension". Note that, by definition, the integral (1.2.2) is 0 if the degree of the integrand does not match the dimension of the virtual fundamental class.

Typically, the cohomology classes $\gamma_1, \ldots, \gamma_k$ are the Poincaré duals of algebraic cycles $\Gamma_1, \ldots, \Gamma_k$ on $X$; in this case, the condition that the degrees of $\operatorname{ev}_1^* \gamma_1, \ldots, \operatorname{ev}_k^* \gamma_k$ sum up to $2\operatorname{vdim} \overline{M}_{g,k}(X, \beta)$ is equivalent to the equality

$$\sum\nolimits_{i=1}^k \bigl(\operatorname{codim}_X(\Gamma_i) - 1\bigr) = \operatorname{expdim} \overline{M}_{g,0}(X, \beta),$$

---

[1] at a point of $\overline{M}_{g,k}(X, \beta)$ corresponding to a map $f : C \to X$ such that $C$ is smooth and $f$ is unramified, this may be computed as $\chi(N_f) + k$, where $N_f$ (the *normal sheaf* of $f$) is the locally free cokernel of the injective map of locally free sheaves $T_C \to f^* T_X$, see [**, Section **], [38, § 3.4.2], or [32, § $1\frac{1}{2}$]; these references also explain how to compute the expected dimension in general.

[2] in fact there are other Gromov–Witten invariants, which are more general intersection numbers against $[\overline{M}_{g,k}(X, \beta)]^{\mathrm{vir}}$ than those in (1.2.2) (we will see some examples in Section 5.3), but our primary interest in these notes resides in the latter.

which means that the incidence conditions imposed by $\Gamma_1, \ldots, \Gamma_k$ to genus $g$ curves in the class $\beta$ are expected to define a finite number of curves. Then, under suitable transversality assumptions, and provided the moduli space $\overline{M}_{g,k}(X, \beta)$ (or equivalently $\overline{M}_{g,0}(X, \beta)$) has the expected dimension, the Gromov–Witten invariant (1.2.2) gives the number of genus $g$ curves in the class $\beta$ (interpreted as stable maps, and counted with multiplicities) which pass through the cycles $\Gamma_1, \ldots, \Gamma_k$. We will be mainly concerned in this text with the case when all $\Gamma_i$'s are points, which is the only relevant case when $X$ is a surface.

**(1.3) Terminology and conventions.** We always work over the field of complex numbers.

Let $C$ be a curve. Its arithmetic genus, denoted by $p_a(C)$, is the integer $1 - \chi(\mathcal{O}_C)$. If $C$ is reduced, its geometric genus is the arithmetic genus of its normalization, and is denoted by $p_g(C)$. When I write 'genus', this means 'geometric genus'.

A reduced curve $C$ is *immersed* when the differential of its normalization map is everywhere non-degenerate. Concretely this means that $C$ has no cuspidal points; it may have however points of any multiplicity, and non-ordinary singularities (e.g., a tacnode, i.e., a point at which there are two smooth local branches tangent one to another). A *node* is an ordinary double point.

A $K3$ surface $S$ is a smooth surface with trivial canonical bundle and vanishing irregularity (i.e., $q(S) = h^1(S, \mathcal{O}_S) = 0$); we may occasionally qualify as $K3$ a surface with canonical singularities, the minimal smooth model of which is a smooth $K3$ surface. Let $p$ be a positive integer. A $K3$ surface of genus $p$ is a pair $(S, L)$, where $S$ is a $K3$ surface and $L$ an effective line bundle on $S$, such that $L^2 = 2p - 2$ (in particular, the $K3$ surface $S$ is algebraic). Under these assumptions, the complete linear system $|L|$ has dimension $p$, and its general member has arithmetic genus $p$ (it is smooth if $|L|$ is free of fixed components). The pair $(S, L)$ is *primitive* if the line bundle $L$ is indivisible, i.e., there is no line bundle $L'$ on $S$ such that $L \cong (L')^{\otimes m}$ for some integer $m > 1$.

In the notation of (1.2), we write $\overline{M}_g(X, \beta)$ for $\overline{M}_{g,0}(X, \beta)$. If $S$ is a surface equipped with an effective line bundle $L \to S$, we write $\overline{M}_{g,k}(X, L)$ for $\overline{M}_{g,k}(X, \beta)$ where $\beta$ is the homology class of the members of $|L|$.

**(1.4)** Let $(S, L)$ be a $K3$ surface of genus $p$. Members of $|L|$ with exactly $\delta$ nodes as singularities have geometric genus $p - \delta$, and are expected to fill up a locus of codimension $\delta$ in $|L|$. For this reason (and because $|L|$ has dimension $p$), the locus of genus $g$ curves in $|L|$ has expected dimension $g$; note that this does not match with the virtual dimension (1.2.1) of $\overline{M}_g(S, L)$, see (3.1). One can actually prove that this is indeed the correct dimension, and that the locus of genus $g$ curves is equidimensional (see [12, § 4.2]). This implies that for a general set of $g$ points $x_1, \ldots, x_g \in S$, there is a finite number of genus $g$ curves in $|L|$ passing through all points $x_1, \ldots, x_g$, by Kleiman and Piene's little lemma [**, Lemma **].

**(1.5) Acknowledgements.** I thank Jim Bryan and Rahul Pandharipande for patiently answering my naive questions, and Yifan Zhao for his observations and help on Section 2.4.

## 2 – Rational curves in a primitive class

In this section we discuss the following result, proved by Beauville [4] following a strategy proposed by Yau and Zaslow [44]. This strategy was inspired by physics; it is an elaboration of the elementary argumentation using Euler numbers presented in Section 2.1.

**(2.1) Theorem** (Yau–Zaslow, Beauville)**.** *Let $(S, L)$ be a smooth primitive $K3$ surface of genus $p_0$, and assume that $\mathrm{Pic}(S) = \mathbf{Z} L$. Then there is a finite number $N^{p_0}$ of rational curves in the complete linear system $|L|$, and it is determined by the formula*

$$(2.1.1) \qquad 1 + \sum_{p=1}^{+\infty} N^p q^p \;=\; \prod_{n=1}^{+\infty} \frac{1}{(1-q^n)^{24}}$$
$$= \; 1 + 24q + 324q^2 + 3200q^3 + \cdots$$

As one may expect, the number $N^{p_0}$ is the number of rational curves counted with appropriate multiplicities. As we shall see, the multiplicity of an integral rational curve $C$ is the topological Euler number $e(\bar{J}C)$ of its compactified Jacobian, and it depends only on the equisingularity class of its singularities. It may be explicitly computed in many cases; in particular, it is 1 whenever the curve is immersed. For a very general $(S,L)$, all rational curves in $|L|$ are nodal by [8], hence in this case all rational curves count with multiplicity one.

The assumption that $\mathrm{Pic}(S) = \mathbf{Z} L$ is made to ensure that all members of $|L|$ are integral curves. While it may be possible to adapt the arguments given by Beauville and Yau–Zaslow to the case when some members of $|L|$ are reducible but reduced, the presence of non-reduced members of $|L|$ poses problems of a different nature, see Section 4.

Theorem (2.1) is a particular case of the more general result discussed in Section 3.2. We will recall there the relevant facts from the theory of modular forms needed to explore the modular aspects of formula (2.1.1), and give more values of $N^p$ for small $p$.

## 2.1 – An elementary topological counting formula

Let $X$ be a complex compact manifold and $\mathfrak{f}$ a 1-dimensional family of divisors of $X$, such that the general member of $\mathfrak{f}$ is smooth and $\mathfrak{f}$ has no fixed part. It is possible to count the number of singular members of $\mathfrak{f}$ using the topological formula stated in Lemma (2.3) below.

**(2.2) Euler number.** Let $X$ be a topological space. The *(topological) Euler number* of $X$ is

$$e(X) = \sum_i (-1)^i \dim H_c^i(X, \mathbf{Z}),$$

where the subscript $c$ means cohomology with compact support.

If $X$ is compact and $F \subseteq X$ is a closed subset, then the following *additivity formula* holds:

$$e(X) = e(X - F) + e(F).$$

It follows from the long exact sequence

$$\cdots \to H_c^i(X - F, \mathbf{Z}) \to H^i(X, \mathbf{Z}) \to H^i(F, \mathbf{Z}) \to H_c^{i+1}(X - F, \mathbf{Z}) \to \cdots$$

see [18, Théorème 4.10.1].

**(2.3) Lemma.** *Let $f : X \to C$ be a surjective morphism from a complex compact manifold onto a smooth curve. One has*

$$(2.3.1) \qquad e(X) = e(F_{\mathrm{gen}})\, e(C) + \sum_{y \in \mathrm{Disc}(f)} \bigl( e(F_y) - e(F_{\mathrm{gen}}) \bigr),$$

*where $F_{\mathrm{gen}}$ is the fibre of $f$ over a general point of $B$, and for all $y \in B$, $F_y$ is the fibre over $y$; $\mathrm{Disc}(f) \subseteq C$ is the set of points above which $f$ is not smooth.*

*Proof.* Set $U = X - \bigcup_{y \in \mathrm{Disc}(f)} F_y$. Then $f : U \to C - \mathrm{Disc}(f)$ is a topological fibre bundle, hence
$$e(U) = e(C - \mathrm{Disc}(f)) \, e(F_{\mathrm{gen}}).$$
Now, by the additivity formula,
$$e(X) = e(U) + \sum_{y \in \mathrm{Disc}(f)} e(F_y) \quad \text{and} \quad e(C - \mathrm{Disc}(f)) = e(C) - \sum_{y \in \mathrm{Disc}(f)} 1,$$
and Formula (2.3.1) follows. □

**(2.4)** The upshot of Lemma (2.3) is that the number of singular fibres of $f : X \to C$ is $e(X) - e(F_{\mathrm{gen}}) e(C)$, provided each singular fibre $F_y$ is counted with multiplicity $e(F_y) - e(F_{\mathrm{gen}})$. When $X$ is a surface and the schematic fibre over $y$ is reduced, the difference $e(F_y) - e(F_{\mathrm{gen}})$ is determined by the singularities of $F_y$, as we shall now explain.

Let $D$ be a reduced projective curve, $\Sigma$ its singular locus, $\nu : \bar{D} \to D$ its normalization, and $\bar{\Sigma} = \nu^{-1}(\Sigma)$. Since $\bar{D} - \bar{\Sigma}$ and $D - \Sigma$ are isomorphic one has, by the additivity formula,

(2.4.1)
$$\begin{aligned} e(D) &= e(D - \Sigma) + e(\Sigma) \\ &= e(\bar{D} - \bar{\Sigma}) + e(\Sigma) \\ &= e(\bar{D} - \bar{\Sigma}) + e(\bar{\Sigma}) + e(\Sigma) - e(\bar{\Sigma}) = e(\bar{D}) - \big(\mathrm{Card}(\bar{\Sigma}) - \mathrm{Card}(\Sigma)\big). \end{aligned}$$

Now, let $f : S \to C$ be a surjective morphism from a smooth projective surface to a smooth curve, and consider a point $y \in C$ such that the schematic fibre $F_y$ is reduced. Then the curve $F_y$ has the same arithmetic genus as the general fibre $F_{\mathrm{gen}}$, and $\bar{F}_y$ is a smooth curve of genus $p_a(F_{\mathrm{gen}}) - \delta$, where $\delta = p_a(F_y) - p_g(F_y)$ is the sum of the $\delta$-invariants of all singularities of $F_y$. Therefore,
$$e(\bar{F}_y) = 2 - 2(p_a(F_{\mathrm{gen}}) - \delta) = e(F_{\mathrm{gen}}) + 2\delta,$$
hence
$$e(F_y) - e(F_{\mathrm{gen}}) = 2\delta - \big(\mathrm{Card}(\bar{\Sigma}_y) - \mathrm{Card}(\Sigma_y)\big)$$
by (2.4.1). This proves the following.

**(2.5) Lemma.** *In the above notation, the multiplicity with which the fibre $F_y$ is counted in Formula* (2.3.1) *is a sum of local multiplicities computed at the singular points of $F_y$, namely*
$$e(F_y) - e(F_{\mathrm{gen}}) = \sum_{x \in \mathrm{Sing}\, F_y} \left( \delta(F_y, x) - \#\binom{\text{local branches}}{\text{of } F_y \text{ at } x} + 1 \right).$$

This gives for example the local multiplicity 1 for a node, 2 for an ordinary cusp, 3 for a tacnode, and 4 for an ordinary triple point.

**(2.6) Application to $K3$ surfaces of genus one.** Let $(S, L)$ be a primitive $K3$ surface of genus one. If $\mathrm{Pic}(S) = \mathbf{Z}\, L$, then $|L|$ is a base-point-free pencil of elliptic curves, and all its members are reduced since $L$ is not divisible. All $K3$ surfaces $S$ have $e(S) = 24$, and smooth elliptic curves have Euler number 0. Thus, it follows from formula (2.3.1) that $|L|$ has 24 singular members, counted with the multiplicity given in Lemma (2.5). This agrees with Theorem (2.1).

**(2.7)** Finally, if $\mathfrak{f}$ is a 1-dimensional family of divisors on a compact complex manifold $X$, free of fixed component and such that the general member of $\mathfrak{f}$ is smooth, in order to find the number of singular members of $\mathfrak{f}$ one may consider the minimal smooth blow-up $\tilde{X} \to X$ such that the proper transform of $\mathfrak{f}$ is base-point-free; then there is a regular map $\tilde{X} \to \mathfrak{f}$, the fibres of which are exactly the members of $\mathfrak{f}$, and one may apply Lemma (2.3) to this fibration.

## 2.2 – Proof of the Beauville–Yau–Zaslow formula

In this section, we outline the proof of Theorem (2.1). Let $p$ be a positive integer, $(S, L)$ a smooth primitive $K3$ surface of genus $p$, and $\mathfrak{L}$ the complete linear system $|L|$. We assume that all members of $\mathfrak{L}$ are integral.

We let $\mathcal{C}$ be the universal curve over $\mathfrak{L}$, and consider the component

$$\pi : \bar{\mathcal{J}}^p \mathcal{C} \to \mathfrak{L}$$

of the compactified Picard scheme of the family $\mathcal{C} \to \mathfrak{L}$ which parametrizes pairs $(C, M)$ made of a member $C$ of $\mathfrak{L}$ and a rank 1, torsion-free, coherent sheaf $M$ of degree $p$ on $C$. For all $[C] \in \mathfrak{L}$, the fibre $\pi^{-1}([C])$ is the compactified Jacobian $\bar{J}^p C$ of $C$. The total space $\bar{\mathcal{J}}^p \mathcal{C}$ is a projective variety of dimension $2p$, which we will call the universal compactified Jacobian.

**(2.8) Claim.** The Euler number $e(\bar{\mathcal{J}}^p \mathcal{C})$ is the number of rational curves in $\mathfrak{L}$, where each rational curve $C$ is counted with multiplicity $e(\bar{J}^p C)$.

*Proof.* The claim relies on the fact that for all $[C] \in \mathfrak{L}$, the fibre $\pi^{-1}([C]) = \bar{J}^p C$ has non-vanishing Euler number if and only if $C$ is rational, see Propositions (2.10) and (2.13) in Section 2.3 below.

The rest is essentially an elaboration of Lemma (2.3) and its proof. There exists a stratification $\mathfrak{L} = \coprod_\alpha \Sigma_\alpha$ by locally closed subsets such that $\pi$ is topologically locally trivial over each stratum $\Sigma_\alpha$ [42]. Then,

$$(2.8.1) \qquad e(\bar{\mathcal{J}}^p \mathcal{C}) = \sum_\alpha e\bigl(\pi^{-1}(\Sigma_\alpha)\bigr)$$

by the additivity formula, and for each $\alpha$ one has

$$e\bigl(\pi^{-1}(\Sigma_\alpha)\bigr) = e(\Sigma_\alpha) \times e(\bar{J}_\alpha)$$

where $\bar{J}_\alpha$ stands for the fibre of $\pi$ over any point in the stratum $\Sigma_\alpha$. Since $e(\bar{J}^p C) = 0$ if the curve $C$ is not rational, the sum in (2.8.1) reduces to the sum over all $\alpha$ such that the curves parametrized by the stratum $\Sigma_\alpha$ are rational. Now, the union of these strata is a finite set (see, e.g., [12, Prop. (4.7)]), which we shall denote by $\mathfrak{L}_{\mathrm{rat}}$, and we finally find

$$(2.8.2) \qquad e(\bar{\mathcal{J}}^p \mathcal{C}) = \sum_{[C] \in \mathfrak{L}_{\mathrm{rat}}} e(\bar{J}^p C),$$

which proves the claim. □

Let $S^{[p]}$ be the component of the Hilbert scheme of $S$ parametrizing 0-dimensional subschemes of length $p$. It is a smooth Hyperkähler variety of dimension $2p$, see [3].

**(2.9) Claim.** The universal compactified Jacobian $\bar{\mathcal{J}}^p \mathcal{C}$ and the Hilbert scheme $S^{[p]}$ have the same Euler numbers, given by Formula (2.1.1).

*Proof.* First, the universal compactified Jacobian $\bar{\mathcal{J}}^p \mathcal{C}$ is smooth and Hyperkähler, as noted in [30, Example 0.5], because it is a connected component of the moduli space of simple sheaves on the $K3$ surface $S$. It is sometimes called the Beauville–Mukai system.

Next, let us prove that $\bar{\mathcal{J}}^p \mathcal{C}$ and $S^{[p]}$ are birational. This will imply, by a theorem of Huybrechts [22, p. 65], that $\bar{\mathcal{J}}^p$ and $S^{[p]}$ are deformation equivalent, hence share the same Betti numbers, so that $e(\bar{\mathcal{J}}^p) = e(S^{[p]})$.

There is an open subset $U \subseteq \bar{\mathcal{J}}^p \mathcal{C}$ whose points are pairs $(C, M)$ such that $C$ a smooth (genus $p$) member of $\mathfrak{L}$ and $M$ is a non-special, degree $p$, line bundle on $C$. For such a pair one has $h^0(C, M) = 1$, and the unique divisor $D$ in the complete linear system $|M|$ may be seen as a point of $S^{[p]}$: this gives a regular map $U \to S^{[p]}$. On the other hand, the fact that $h^1(C, \mathcal{O}_C(D)) = 0$ implies by geometric Riemann–Roch that $D$ spans the whole space in the canonical embedding of $C$, hence $C$ is the only member of $\mathfrak{L}$ that contains $D$. Therefore, the map $(C, M) \in U \mapsto D$ is injective, and this ends the proof of the birationality of $\bar{\mathcal{J}}^p \mathcal{C}$ and $S^{[p]}$.

Finally, the fact that $e(S^{[p]})$ is the coefficient of $q^p$ in the Fourier expansion of $\prod_{n=1}^{\infty}(1-q^n)^{-24}$ as in (2.1.1) follows from the computation by Göttsche of the Betti and Euler numbers of $S^{[p]}$ for any smooth, complex, projective surface [20]. This computation is based on the idea of using the Weil conjectures (proved by Deligne) to translate the question into that of counting the points of $S^{[p]}$ over finite fields. □

Claims (2.8) and (2.9) together imply Theorem (2.1).

## 2.3 – Compactified Jacobian of an integral curve

Let $C$ be an integral curve. We now turn to the study of the Euler number of the compactified Jacobian $\bar{J}^d C$. We will in particular prove some results used in the proof of Claim (2.8) in Section 2.2. As is well-known, the choice of an invertible sheaf of degree $d$ on $C$ induces an isomorphism between $\bar{J}^d C$ and $\bar{J}C = \bar{J}^0 C$, so we will actually consider the latter variety.

**(2.10) Proposition.** *If $C$ is an integral curve of positive geometric genus, then $e(\bar{J}C) = 0$.*

*Proof.* There is an exact sequence

$$0 \to H \to JC \to J\tilde{C} \to 0$$

where $\tilde{C}$ is the normalization of $C$,[3] $J$ denotes the Jacobian $\text{Pic}^0$, and $H$ is a product of copies of $(\mathbf{C}, +)$ and $(\mathbf{C}^*, \times)$ (see, e.g., [25, Thm. 7.5.19]). Since $H$ is divisible, the sequence splits as an exact sequence of Abelian groups, hence $JC = J\tilde{C} \times H$. Since $J\tilde{C}$ is a positive dimensional torus, for all positive integers $n$, there exists an order $n$ group $G_n$ which is a subgroup of both $JC$ and $J\tilde{C}$.

Then, [4, Lemma 2.1] tells us that for all $M \in G_n$ and $\mathcal{F} \in \bar{J}C$, the two sheaves $\mathcal{F}$ and $\mathcal{F} \otimes M$ are not isomorphic. Thus $G_n$ acts freely on $\bar{J}C$, hence $n$ divides $e(\bar{J}C)$. Since this is true for all $n$, we conclude that $e(\bar{J}C) = 0$. □

We need the following local construction (cf. [4, §3.6] and the references therein) in order to make $e(\bar{J}C)$ explicit for a rational curve $C$.

**(2.11)** Let $(C, x)$ be a germ of curve which we assume to be unibranch (see (2.12) below), and let $\tilde{C}$ be the normalization of $C$; there is only one point in the preimage of $x$, which we call $x$ as well. Set $\delta_x = \dim \mathcal{O}_{\tilde{C},x} / \mathcal{O}_{C,x}$ (this is the number by which a singularity equivalent to $(C, x)$ makes the geometric genus drop), and let $\mathbf{c}'_x$ be the ideal $\mathcal{O}_{\tilde{C}}(-2\delta \cdot x)$.[4] We shall consider the two finite-dimensional algebras $A_x := \mathcal{O}_{C,x}/\mathbf{c}'$ and $\tilde{A}_x := \mathcal{O}_{\tilde{C},x}/\mathbf{c}'$.

Let $\mathbf{G}_x$ be the closed subvariety of the Grassmannian $\mathbf{G}(\delta_x, \tilde{A}_x)$ parametrizing those codimension $\delta_x$ subspaces of $\tilde{A}_x$ that are also sub-$A_x$-modules of $\tilde{A}_x$; $\mathbf{G}_x$ may also be seen as the variety parametrizing codimension $\delta_x$ sub-$\mathcal{O}_{C,x}$-modules of $\mathcal{O}_{\tilde{C},x}$. It only depends on the completion $\hat{\mathcal{O}}_{C,x}$, hence only on the analytic type of the singularity $(C, x)$.

---

[3] exceptionally, I do not use the notation $\bar{C}$ in order to avoid unpleasant confusions between $\bar{J}C$ and $J\bar{C}$

[4] we reserve the notation $\mathbf{c}_x$ for the conductor ideal (see [**]), which contains $\mathbf{c}'_x$

**(2.12)** A curve is *unibranch* if its normalization is a homeomorphism, or equivalently if it is everywhere analytically locally irreducible. All curves $C$ have a *unibranchization* $\check{\nu}: \check{C} \to C$: this is the unique partial normalization of $C$ such that $\check{C}$ is unibranch, and such that all others partial normalizations $\nu': C' \to C$ with $C'$ unibranch factor through $\check{\nu}$.

If $C$ is a unibranch curve with singular locus $\Sigma \subseteq C$, the product $\prod_{x \in \Sigma} \mathbf{G}_x$ parametrizes sub-$\mathcal{O}_C$-modules $\mathcal{F} \subseteq \mathcal{O}_{\tilde{C}}$ such that $\dim \mathcal{O}_{\tilde{C},x}/\mathcal{F}_x = \delta_x$ for all $x$. Such an $\mathcal{F}$ enjoys the property that $\chi(\mathcal{F}) = \chi(\mathcal{O}_C)$, which implies $\mathcal{F} \in \bar{J}C$. This defines a morphism

$$\varepsilon: \prod_{x \in \Sigma} \mathbf{G}_x \to \bar{J}C.$$

Part (ii) in the following statement is proved by showing that $\varepsilon$ is a homeomorphism if $C$ is rational, though in general not an isomorphism.

**(2.13) Proposition.** (i) [4, Prop. 3.3] *If $C$ is an integral curve, then $e(\bar{J}C) = e(\bar{J}\check{C})$.*
(ii) [4, Prop. 3.8] *If $C$ is a unibranch rational curve with singular locus $\Sigma$, then $e(\bar{J}C) = \prod_{x \in \Sigma} e(\mathbf{G}_x)$.*

By (i), one has $e(\bar{J}C) = 1$ for an immersed rational curve $C$ ($C$ is immersed if all its analytic local branches are smooth, equivalently if the unibranchization $\check{C}$ is smooth). If $C$ is not immersed, then $e(\bar{J}C) = e(\bar{J}\check{C})$ and the latter number is computed by (ii); it follows from Theorem (2.17) in the next subsection that this number is positive if $C$ is a rational curve lying on a smooth surface.

Part (ii) also shows that for all rational curves $C$ (unibranch or not, thanks to (i)), $e(\bar{J}C)$ is entirely determined by the singularities of $C$. In fact it only depends on their topological types if $C$ has planar singularities, see Corollary (2.18) below.

**(2.14) Example** ([4, § 4]). If $(C,x)$ is the germ of curve given by the equation $u^p + v^q = 0$ at the origin in the affine plane, with $p$ and $q$ relatively prime, then

$$(2.14.1) \qquad e(\mathbf{G}_x) = \frac{1}{p+q}\binom{p+q}{p}.$$

Formula (2.14.1) also gives $e(\mathbf{G}_{x'})$ for all singularities $(C',x')$ equisingular to $(C,x)$ by the aforementioned Corollary (2.18).

In particular, $e(\mathbf{G}_x) = \ell + 1$ for the cuspidal singularity defined by the equation $u^2 + v^{2\ell+1} = 0$, i.e., for an $A_{2\ell}$-double point; on the other hand $A_{2\ell+1}$-double points are immersed, hence $e(\mathbf{G}_x) = 1$ for them. Using Proposition (2.13), (i), to the effect that the local contribution $e(\mathbf{G}_x)$ of a germ $(C,x)$ is the product of the local contributions of all analytic local irreducible branches of $(C,x)$, Formula (2.14.1) is enough to determine the local contribution $e(\mathbf{G}_x)$ of any simple curve singularity, see [4, Prop. 4.5].

**(2.15) Remark.** The fact that all immersed rational curves count with multiplicity one seems to disagree with the results of Section 2.1, see in particular Lemma (2.5). However, if $|F|$ is a complete pencil of elliptic curves, the assumption that all curves in $|F|$ are integral implies that all rational curves in $|F|$ are curves of arithmetic genus 1 with either a node or an ordinary cusp as their unique singular point, and in this case the two multiplicities predicted by Proposition (2.13) and Lemma (2.5) do agree.

**(2.16) Remark.** Returning to the case of a $p$-dimensional linear system $\mathfrak{L}$ with all members integral, Beauville deems "rather surprising" the fact that "some highly singular [immersed] curves count with multiplicity one", and considers the case $p = 2$ to provide a confirming example.

*(2.16.1)* In hindsight, a good conceptual explanation comes from the fact that the numbers $N^p$ should be seen as counting stable maps rather than embedded curves, since stable maps don't make any difference between nodal and arbitrary immersed curves: if $p$ is a singular point of a curve $C$ and $p'$ is a point of the normalization $\tilde{C}$ mapping to $p$, analytically locally at $p'$ one does not see how the local branch of $C$ at $p$ corresponding to $p'$ intersects the other local branches (one does see however if this branch of $C$ is singular, by looking at the differential of the map $\tilde{C} \to C$). This heuristic argument is confirmed by Theorem (2.20) in the next subsection, see also Lemma (3.4).

*(2.16.2)* Another comforting observation is that the equigeneric deformation space of an immersed singularity is smooth, as we shall see in (2.19) in the next subsection. Indeed, let $C$ be a rational member of $\mathfrak{L}$. There exists an analytic neighbourhood of $[C]$ in $\mathfrak{L}$ and a map $\phi$ from $U$ to the semi-universal deformation space $B$ of the singularities of $C$, which encodes the singularities of the curves neighbouring $C$ in $\mathfrak{L}$. Then, the locus of rational curves in $U$ is the inverse image by $\phi$ of the equigeneric locus $\mathrm{EG}(C)$ in $B$, hence it is reasonable to imagine that the multiplicity of $C$ in the count of rational curves is the multiplicity of the 0-cycle $\phi^*\mathrm{EG}(C)$ at the point $[C] \in U$. The fact that rational curves in $\mathfrak{L}$ are finitely many implies as in (1.4) that the intersection $\phi(U) \cap \mathrm{EG}(C)$ has the expected dimension. If we assume furthermore that $\phi$ is an immersion and $\phi(U)$ intersects $\mathrm{EG}(C)$ transversely, then the multiplicity of $\phi^*\mathrm{EG}(C)$ at $[C]$ equals the multiplicity of the origin of $B$ as a point of $\mathrm{EG}(C)$.

In the upshot, it is plausible that a rational curve counts with multiplicity equal to the multiplicity of the origin in the equigeneric deformation space of its singularities. Since the latter multiplicity is one if all singularities of $C$ are immersed, it is plausible that any immersed rational curve counts with multiplicity one. In fact, we shall see in Section 2.4 that things happen exactly as described above.

## 2.4 – Two fundamental interpretations of the multiplicity

We now conclude Section 2 by stating two enlightening geometric interpretations of the local multiplicities $e(\mathbf{G}_x)$ defined in the previous subsection, due to Fantechi, Göttsche and van Straten [14].

In the present subsection we only consider *planar* curve singularities $(C, x)$, i.e., singularities that may be realized by the germ of a plane curve. All singularities occuring on a curve that lies on a smooth surface are of this kind.

**(2.17) Theorem** ([14, Thm. 1])**.** *Let $(C, x)$ be a reduced, planar, curve singularity, and let $\mathbf{G}_x$ be as in (2.11). Then the topological Euler number $e(\mathbf{G}_x)$ equals the multiplicity at the point $[(C, x)]$ of the equigeneric locus $\mathrm{EG}(C, x)$ in the semi-universal deformation space of the singularity $(C, x)$.*

Recall that the *equigeneric locus* $\mathrm{EG}(C, x)$ is defined as the reduced subscheme of the semi-universal deformation space of the singularity $(C, x)$ supported on those points corresponding to singularities with the same $\delta$-invariant as $(C, x)$. In practice this is the closure of the locus along which the singularity deforms to $\delta$ nodes, where $\delta$ is the number by which the singularity makes the geometric genus drop. See, e.g., [**] for more details.

Together with Proposition (2.13), Theorem (2.17) implies that for a rational, integral, curve $C$ with only planar singularities, the Euler number $e(\bar{J}C)$ of the compactified Jacobian of $C$ equals the multiplicity at the point $[C]$ of the equigeneric stratum $\mathrm{EG}(C)$ in a semi-universal deformation space of the singularities of $C$. See (2.16.2) above for a discussion of why this provides a reasonable notion of multiplicity for $C$ in a count of rational curves.

Theorem (2.17) has been generalized by Shende [40], to the effect that the multiplicity at $[C]$ of all $\delta$-constant strata in a semi-universal deformation space of the singularities of $C$ are determined by the Euler numbers of the Hilbert schemes of points on $C$ (for a curve $C$ of arbitrary genus).

On the other hand, Maulik [27] has proved that the Euler numbers of the Hilbert schemes of points $(C,x)^{[n]}$ of a planar curve singularity only depend on the topological type of the singularity. This implies the following.

**(2.18) Corollary** (see [27, § 6.5]). *Let $(C,x)$ be a reduced, planar, curve singularity. The multiplicities at the point $[(C,x)]$ of the $\delta$-constant (a.k.a. constant genus) strata of the semi-universal deformation space of the singularity $(C,x)$ only depend on the topological type (or equisingularity class) of $(C,x)$.*

**(2.19)** As a corollary of Theorem (2.17) and Proposition (2.13), one obtains that if $(C,x)$ is an immersed planar curve singularity, then the equigeneric locus $\mathrm{EG}(C,x)$ in the space of semi-universal deformations is smooth at the point $[(C,x)]$. Indeed, note that for all reduced plane curve singularities $(C,x)$, there exists a rational curve $\tilde{C}$ with $(C,x)$ as its only singularity (this follows for instance from [29]), and then one has $e(\mathbf{G}_x) = e(\bar{J}\tilde{C})$.

In fact the converse holds as well, i.e., for a reduced planar curve singularity $(C,x)$, the equigeneric locus $\mathrm{EG}(C,x)$ is smooth at the point $[(C,x)]$ if and only if $(C,x)$ is immersed, see [13, Proposition (4.17), (2)].

Now let us turn to the second interpretation of the multiplicity $e(\mathbf{G}_x)$, or rather of $e(\bar{J}C) = \prod_{x \in \Sigma} e(\mathbf{G}_x)$. It tells us that what the Yau–Zaslow formula (2.1.1) really computes is the numbers of genus 0 stable maps realizing primitive classes on $K3$ surfaces. Recall that if $C$ is an integral curve of geometric genus $g$ which is a subvariety of a complex variety $X$, then $\overline{M}_g(X,[C])$ is the space of genus $g$ stable maps with target $X$ realizing the homology class $[C] \in H_2(C,\mathbf{Z})$.

**(2.20) Theorem** ([14, Thm. 2] and [45]). *Let $C$ be an integral rational curve contained in a smooth $K3$ surface $S$. Then the topological Euler number $e(\bar{J}C)$ equals the length of the space of stable maps $\overline{M}_0(S,[C])$ at the closed point corresponding to the normalization of $C$.*

I have been informed by Yifan Zhao that he has found a problem in the proof of [14, Thm. 2], and will provide a corrected proof in the upcoming [45]; I refer to [45] for the full details.

Fantechi–Göttsche–van Straten first consider the space $\overline{M}_g(C,[C])$ of stable maps to $C$ itself which realize the fundamental class $[C]$. It is a 0-dimensional scheme, which contains a single closed point, corresponding to the normalization $\nu : \bar{C} \to C$. If $C$ is an isolated, genus $g$ curve in a smooth manifold $X$, then $\overline{M}_g(C,[C])$ is a subscheme of $\overline{M}(X,[C])$ and, at the normalization of $C$, the length of the former scheme is a lower bound for the length of the latter.

Now, in [14, Thm. 2] the authors claim that for an integral projective curve $C$ of genus $g$ with planar singularities, the length of $\overline{M}_g(C,[C])$ equals the multiplicity at $[C]$ of the semi-universal deformation space of the singularities of $C$. However, [31, §7.2.2] provides a counterexample, namely a rational curve $C$ with planar singularities such that the Euler number $e(\bar{J}C)$ is larger than the length of $\overline{M}_0(C,[C])$ (as we have seen, since $C$ is rational, $e(\bar{J}C)$ equals the multiplicity at $[C]$ of the semi-universal deformation space).

# 3 – Curves of any genus in a primitive class

## 3.1 – Reduced Gromov–Witten theories for $K3$ surfaces

**(3.1) A vanishing phenomenon.** It happens that all Gromov–Witten invariants of $K3$ surfaces are trivial. The fundamental reason for this is that Gromov–Witten invariants are defor-mation invariant (this is indeed a desirable feature of any well-behaved counting invariants, and the very reason why they are called invariants and not merely numbers, or integrals), and all algebraic $K3$ surfaces deform to non-algebraic ones, which in general do not contain any curve at all.

Somewhat more concretely, the explanation is that the virtual and actual dimensions of the moduli spaces of stable maps on $K3$ surfaces do not match, as we already pointed out in (1.4). Let $S$ be a $K3$ surface, and $C \subseteq S$ an integral curve of geometric genus $g$. Consider the stable map $f : \bar{C} \to S$ obtained by composing the normalization $\bar{C} \to C$ with the inclusion $C \subseteq S$. The normal sheaf $N_f$ of $f$ is isomorphic to the canonical bundle $\omega_{\bar{C}}$, and therefore $h^0(N_f) = g$ and $h^1(N_f) = 1$. It follows that the virtual dimension of $\overline{M}_g(S, [C])$ is $g - 1$, whereas the curve $C$ actually moves in a $g$-dimensional family of curves of genus $g$ on $S$ (see [12, § 4.2] for more details). This implies that the Gromov–Witten invariant counting genus $g$ curves on $S$ passing through the appropriate number of points (namely $g$) vanishes for mere degree reasons.

The following two paths have been successfully followed to circumvent this issue, and define modified invariants for algebraic $K3$ surfaces which capture the relevant enumerative information.

**(3.2) Invariants of families of symplectic structures** [6, § 2–3]**.** This was chronologically the first workaround to be proposed, and it enabled the counting of curves of any genus in a primitive class on a $K3$ surface reported on in Section 3.2 below.

Let $S$ be a polarized $K3$ surface. The idea here is really to take into account the existence of non-algebraic deformations of $S$. To this effect, instead of counting curves directly on $S$, one counts curves in a family of Kähler surfaces defined over the 2-sphere $\mathbf{S}^2$, canonically attached to $S$, and in which roughly speaking $S$ is the only one to be algebraic, so that all the curves we count are actually concentrated on $S$.

This family of Kähler surfaces is the twistor family of $S$ (cf. [37, p. 124]): the polarization on $S$ determines a Kähler class $\alpha$, and Yau's celebrated theorem asserts that there is a unique Kählerian metric $g$ in $\alpha$ with vanishing Ricci curvature. Then the holonomy defines an action of $\mathbf{H}$ (the field of quaternion numbers) on the holomorphic tangent bundle $TS$ by parallel endomorphisms. The quaternions of square $-1$ define those complex structures on the differentiable manifold $S$ for which the metric $g$ remains Kählerian. There is a 2-sphere worth of such quaternions, and it parametrizes the family we are interested in.

**(3.3) Reduced Gromov–Witten theory** ([28, 2.2], see also [26])**.** In this second approach, the idea is to plug in the fact that, for a stable map $f$ as in (3.1), the space $H^1(\bar{C}, N_f)$ although non-trivial does not contain any actual obstruction to deform $f$ as a map with target $S$, following for instance Ran's results on deformation theory and the semiregularity map. To this end, Maulik and Pandharipande define a suitable perfect obstruction theory which they dub *reduced*, and which provides, following the construction pioneered by Behrend and Fantechi, a reduced virtual fundamental class $[\overline{M}_g(S, \beta)]^{\mathrm{red}}$, for all integers $g \geqslant 0$ and algebraic classes $\beta \in H_2(S, \mathbf{Z})$, which has the appropriate (real) dimension $2g$.

This in turn gives reduced Gromov–Witten integrals, by replacing the virtual fundamental class by its reduced version in the integral (1.2.2). They are invariant under *algebraic* deformations of $K3$ surfaces.

## 3.2 – The Göttsche–Bryan–Leung formula

In this subsection, I discuss a result of Bryan and Leung [6] giving the number of curves in a primitive linear system on a $K3$ surface that have a given genus and pass through the appropriate number of base points. The formula was conjectured by Göttsche [21] as a particular case of a more general framework, see [**].

Let
$$N_g^p := \int_{[\overline{M}_{g,g}(S,L)]^{\mathrm{red}}} \mathrm{ev}_1^*(\mathrm{pt}) \cup \cdots \cup \mathrm{ev}_g^*(\mathrm{pt})$$

be the reduced Gromov–Witten invariant counting curves of genus $g$ in the linear system $|L|$ on a primitive $K3$ surface $(S, L)$ of genus $p$ (pt $\in H^4(S, \mathbf{Z})$ is the point class, and the $\mathrm{ev}_1, \ldots, \mathrm{ev}_g$ are the evaluation maps $\overline{M}_{g,g}(S, L) \to S$, as in (1.2)). The following result tells us that these invariants do indeed count curves.

**(3.4) Lemma.** *The invariants $N_g^p$ are strongly enumerative, in the following sense: let $(S, L)$ be a very general primitive $K3$ surface of genus $p$; then $N_g^p$ is the actual number of genus $g$ curves in $|L|$ passing through a general set of $g$ points, all counted with multiplicity $1$.*

*Proof.* The key fact is that genus $g$ curves on a $K3$ surface move in $g$-dimensional families, see e.g., [12, Prop. (4.7)]; for $g > 0$, this implies by a deformation argument that the general member of such a family is an immersed curve [12, Prop. (4.8)]. The same holds for $g = 0$ by the more difficult result of Chen [8], which asserts that all rational curves in $|L|$ are nodal (this requires the very generality assumption).

Let $\mathbf{x} = (x_1, \ldots, x_g)$ be a general (ordered) set of $g$ points on $S$. The invariant $N_g^p$ may be computed by integration against a virtual class on the so-called cut-down moduli space $\overline{M}(S, \mathbf{x}) \subseteq \overline{M}_{g,g}(S, L)$, consisting of those genus $g$ stable maps sending the $i$-th marked point to $x_i \in S$ for $i = 1, \ldots, g$, see [6, Appendix A].

Let $[f : C \to S] \in \overline{M}(S, \mathbf{x})$. Thanks to the generality assumption on $(S, L)$, we may and will assume that $L$ generates the Picard group of $S$. The condition $f_*C \in |L|$ thus imposes that $f_*C$ is an integral cycle, hence that $f$ contracts all irreducible components of $C$ but one, and restricts to a birational map on the latter component, which we will call $C_1$. Now the points $x_1, \ldots, x_g$, being general, impose $g$ general independent linear conditions on $|L|$, which implies that the curve $f(C) = f(C_1)$ must have geometric genus at least $g$ by the key fact mentioned at the beginning of the proof. Therefore $C_1$ as well must have geometric genus at least $g$, and because of the inequality of arithmetic genera
$$p_a(C_1) \leqslant p_a(C) = g,$$
this implies that $C_1$ is smooth of genus $g$; moreover, the stability conditions then imply that $C = C_1$, hence $f$ is the normalization of the integral genus $g$ curve $f(C_1)$.

This already tells us that the space $\overline{M}(S, \mathbf{x})$ is 0-dimensional, and isomorphic as a set to the space of genus $g$ curves in $|L|$ passing through $x_1, \ldots, x_g$. Since $\overline{M}(S, \mathbf{x})$ has the expected dimension, the virtual class on it simply encodes its schematic structure. It follows that the number $N_g^p$ is weakly enumerative, i.e., it gives the number of genus $g$ curves passing through $x_1, \ldots, x_g$ counted with multiplicities.

The fact that these multiplicities all equal $1$ follows from the fact that all the curves we consider are immersed, as follows. Let $[f : C \to S] \in \overline{M}(S, \mathbf{x})$ as above. Since $f$ is an immersion, it has normal bundle
$$N_f = f^*\omega_S \otimes \omega_C = \omega_C$$
by [12, (2.3)], hence $\mathrm{h}^0(N_f) = g$ and the full moduli space $\overline{M}_{g,0}(S, L)$ is smooth of dimension $g$ at the point $[f]$. By [12, Lemma (2.5)], there is a surjective map $e$ from a neighbourhood of $[f]$ in $\overline{M}_{g,0}(S, L)$ onto a neighbourhood of $f(C)$ in the locally closed subset of $|L|$ parametrizing genus $g$ curves. By generality of the points $x_1, \ldots, x_g$, the latter space is smooth of dimension $g$ at the point $[f(C)]$. Therefore $e$ is a local isomorphism at $[f]$, and this implies that the scheme $\overline{M}(S, \mathbf{x})$ is reduced at $[f]$, which shows that the stable map $f$ counts with multiplicity $1$. □

For all integers $g \geqslant 0$, set

(3.5.1) $$F_g(q) := \sum_{p=g}^{+\infty} N_g^p \, q^p$$

as a formal power series in the variable $q$, where we set $N_0^0$ by convention so that $F_0$ equals the power series of (2.1.1). Beware the shift in degree between the definition (3.5.1) of the series $F_g$ and that given in [6]. Note that $N_g^p = 0$ whenever $p < g$.

**(3.5) Theorem** (Bryan–Leung). *The power series $F_g$ is the Fourier expansion of*

(3.5.2) $$\left(\sum_{n=1}^{+\infty} n\sigma_1(n) \, q^n\right)^g \prod_{m=1}^{+\infty} \frac{1}{(1-q^m)^{24}},$$

*where $\sigma_1(n) := \sum_{d|n} d$ is the sum of all positive integer divisors of $n$.*

This formula gives the possibility to explicitly compute as many numbers $N_g^p$ as we want. Table 1 (p. 14) gives sample values for small $p$ and $g$. Note that since columns are indexed by $\delta := p - g$, the Fourier coefficients of a given $F_g$ are read along a diagonal; this gives for instance $F_0$ as in (2.1.1),

$$F_1(q) = q + 30q^2 + 480q^3 + 5460q^4 + \cdots$$

and so on.

We will discuss the proof of Theorem (3.5) in Section 3.3.

**(3.6) Modularity.** Formula (3.5.2) implies that $qF_g(q)$ is a quasi-modular form for all $g$. In the present paragraph I explain what this means; this is also an occasion to introduce notation for further use (I follow [39]). I will however not try to give a meaningful interpretation of the (quasi-)modularity of $qF_g(q)$.

*(3.6.1)* For every integer $k > 1$, define the $k$-th *Eisenstein series* to be

$$G_k(z) := \sum_{\substack{(m,n) \in \mathbf{Z}^2: \\ (m,n) \neq (0,0)}} \frac{1}{(m+nz)^{2k}};$$

it is a modular form of weight $2k$ [39, Prop. VII.4], which means that it is holomorphic and $G_k(z)dz^k$ is invariant under the action of $\mathrm{PSL}_2(\mathbf{Z})$. Its Fourier expansion at infinity is

$$G_k(z) = \frac{2^{2k}}{(2k)!} B_k \pi^{2k} + 2 \frac{(2\pi i)^{2k}}{(2k-1)!} \sum_{n=1}^{+\infty} \sigma_{2k-1}(n) q^n,$$

where $q = e^{2\pi i z}$, $\sigma_k(n) = \sum_{d|n} d^k$, and $B_k$ is the $k$-th Bernoulli number, defined by the formula

$$\frac{x}{e^x - 1} = 1 - \frac{x}{2} + \sum_{k=1}^{+\infty} (-1)^{k+1} B_k \frac{x^{2k}}{(2k)!}$$

[39, Prop. VII.8]. We set

$$E_k(z) := \frac{(2k)!}{2^{2k} B_k \pi^{2k}} G_k(z) = 1 + (-1)^k \frac{4k}{B_k} \sum_{n=1}^{+\infty} \sigma_{2k-1}(n) q^n$$

$$= 1 + (-1)^k \frac{4k}{B_k} \sum_{n=1}^{+\infty} n^{2k-1} \frac{q^n}{1-q^n}.$$

| $p$ \ $\delta$ | 1 | 2 | 3 | 4 | 5 | 6 | 7 | 8 | 9 |
|---|---|---|---|---|---|---|---|---|---|
| 1 | 24 | | | | | | | | |
| 2 | 30 | 324 | | | | | | | |
| 3 | 36 | 480 | 3200 | | | | | | |
| 4 | 42 | 672 | 5460 | 25650 | | | | | |
| 5 | 48 | 900 | 8728 | 49440 | 176256 | | | | |
| 6 | 54 | 1164 | 13220 | 88830 | 378420 | 1073720 | | | |
| 7 | 60 | 1464 | 19152 | 150300 | 754992 | 2540160 | 5930496 | | |
| 8 | 66 | 1800 | 26740 | 241626 | 1412676 | 5573456 | 15326880 | 30178575 | |
| 9 | 72 | 2172 | 36200 | 371880 | 2499648 | 11436560 | 36693360 | 84602400 | 143184000 |
| 10 | 78 | 2580 | 47748 | 551430 | 4213332 | 22116456 | 81993600 | 219548277 | 432841110 |
| 11 | | 3024 | 61600 | 791940 | 6808176 | 40588544 | 172237344 | 531065070 | 1210781880 |
| 12 | | | 77972 | 1106370 | 10603428 | 71127680 | 342358560 | 1205336715 | 3154067950 |
| 13 | | | | 1508976 | 15990912 | 119665872 | 647773200 | 2582847180 | 7698660544 |
| 14 | | | | | 23442804 | 194196632 | 1172896512 | 5255204625 | 17710394230 |
| 15 | | | | | | 305225984 | 2041899840 | 10205262330 | 38607114200 |
| 16 | | | | | | | 3431986848 | 19002853575 | 80149394030 |
| 17 | | | | | | | | 34070137272 | 159184435520 |
| 18 | | | | | | | | | 303705014550 |

Table 1: First values of $N_g^p$ ($\delta = p - g$)

*(3.6.2)* Define
$$\Delta(z) := \bigl(60G_2(z)\bigr)^3 - 27\bigl(140G_3(z)\bigr)^2,$$
the discriminant of the cubic polynomial $4X^3 - 60G_2X - 140G_3$ divided by 16. It is a modular form of weight 12 vanishing at infinity, and it is a theorem of Jacobi [39, Thm. VII.6] that
$$\Delta(z) = (2\pi)^{12} q \prod_{n=1}^{+\infty} (1-q^n)^{24}.$$

*(3.6.3)* In the $k=1$ case, we set
$$G_1(z) := \sum_{n \in \mathbf{Z}} \sum_{\substack{m \in \mathbf{Z}: \\ (m,n) \neq (0,0)}} \frac{1}{m+nz^2}$$
(note that the order of summation is significant). It has Fourier expansion at infinity
$$G_1(z) = \frac{\pi^2}{3} - 8\pi^2 \sum_{n=1}^{+\infty} \sigma_1(n) q^n$$
[39, § VII.4.4], and we set
$$E_1(z) := \frac{3}{\pi^2} G_1(z) = 1 - 24 \sum_{n=1}^{+\infty} \sigma_1(n) q^n.$$
One has the identity [39, § VII.4.4]
$$\frac{d\Delta}{\Delta} = 2\pi i E_1(z) dz.$$

*(3.6.4)* The function $G_1$ is *not* a modular form, but still it does satisfy a functional equation close to that equivalent to the invariance of $G_1(z)dz$ under the action of $\mathrm{PSL}_2(\mathbf{Z})$ [5, Prop. 6 p. 19]. For this reason, it is called a *quasi-modular* form.

One may then define the ring of quasi-modular forms as the $\mathbf{C}$-algebra generated by $G_1$ and the algebra of modular forms (see [5] for a more intrinsic definition). Since the ring of modular forms is generated by $G_2$ and $G_3$ [39, Cor. VII.2], the ring of quasi-modular forms may be concretely described as $\mathbf{C}[G_1, G_2, G_3]$. The ring of quasi-modular forms is closed under differentiation by the operator
$$D := q\frac{d}{dq} = \frac{1}{2\pi i}\frac{d}{dz}$$
[5, Prop. 15 p. 49].

*(3.6.5) Quasi-modularity of $qF_g(q)$.* Taking into account the various stunning formulæ above, (3.5.2) may be rewritten as
$$qF_g(q) = \left(-\frac{1}{24} q \frac{dE_1}{dq}\right)^g \left(\frac{\Delta}{(2\pi)^{12}}\right)^{-1},$$
from which it follows that $qF_g(q)$ is a quasi-modular form, with a simple pole at infinity (i.e., at $q=0$) if $g=0$.

## 3.3 – Proof of the Göttsche–Bryan–Leung formula

This subsection is dedicated to the proof of Theorem (3.5). As the latter really is about counting actual curves (i.e., subvarieties of the surface), as Lemma (3.4) attests, one may prefer in the first place to avoid the complications of Gromov–Witten theory to prove it. Bryan and Leung's proof however goes the orthogonal direction: it fundamentally relies on the agile possibilities featured in Gromov–Witten theory, precisely as a reward to the aforementioned complications.

**(3.7) Degeneration to an elliptic** $K3$. Let $g \geqslant 0$, $p > 0$ be integers. By deformation invariance, we are free to compute the number $N_g^p$ on our favourite primitively polarized $K3$ surface $(S, L)$ of genus $p$. We let $S$ be an elliptic $K3$ surface with a section $E$, and denote by $F$ the class of the ellitic fibres; the intersection form (or its restriction to the subspace $\langle E, F \rangle$) is given in the basis $(E, F)$ by the matrix

$$\begin{pmatrix} -2 & 1 \\ 1 & 0 \end{pmatrix}.$$

We set $L := \mathcal{O}_S(E + pF)$. Then $L^2 = 2p - 2$, and $L$ is a primitive polarization of genus $p$. We shall compute the numbers $N_g^p$ on the pair $(S, L)$.

**(3.8)** The linear system $|E + pF|$ has dimension $p$ and consists solely of reducible curves $E + F_1 + \cdots + F_p$, where the $F_i$'s are in the class $F$. From this it readily follows that if we fix a general set of $g$ points $\mathbf{y} = (y_1, \ldots, y_g)$ on $S$, then the moduli space

$$\overline{M}_{g,\mathbf{y}}(S, E + pF) := \overline{M}_{g,g}(S, E + pF) \cap \mathrm{ev}_1^*(y_1) \cap \ldots \cap \mathrm{ev}_g^*(y_g)$$

of genus $g$ stable maps passing through the points $y_1, \ldots, y_g$ decomposes as the disjoint union

$$\coprod_{\mathbf{a},\mathbf{b}} \overline{M}_{\mathbf{a},\mathbf{b}},$$

where $\mathbf{a} = (a_1, \ldots, a_{24})$ and $\mathbf{b} = (b_1, \ldots, b_g)$ range through $\mathbf{Z}_{\geqslant 0}^{24}$ and $\mathbf{Z}_{>0}^g$ respectively, subject to the condition that $\sum_i a_i + \sum_j b_j = p$, and $\overline{M}_{\mathbf{a},\mathbf{b}}$ is the moduli space of genus $g$ stable maps $(f : C \to S, x_1, \ldots, x_g)$ such that $f(x_j) = y_j$ for $j = 1, \ldots, g$ and

$$f_*C = E + \sum_{i=1}^{24} a_i R_i + \sum_{j=1}^{g} b_j F_j,$$

where $R_1, \ldots, R_{24}$ are the 24 rational members of the pencil $|F|$, and $F_j$ being the unique member of $|F|$ containing the point $y_j$ for all $j = 1, \ldots, g$ (see figure below); note that we may and do assume that all members of $|F|$ are irreducible, and the $R_i$'s are 1-nodal.

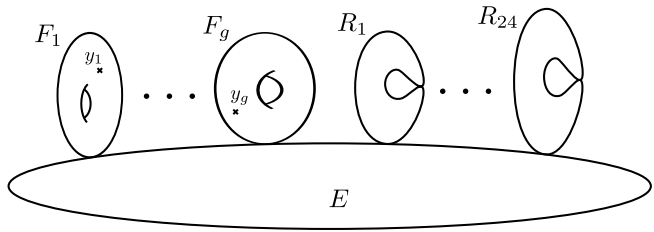

**(3.9) Partition function.** For all positive integers $n$, we let $p(n)$ be the *number of partitions* of $n$, i.e., the number of ways to write $n = \lambda_1 + \cdots + \lambda_k$ with $\lambda_1 \geqslant \cdots \geqslant \lambda_k \geqslant 1$ ($k$ is not fixed). The numbers $p(n)$ may be computed using the generating series

$$1 + \sum_{n=1}^{+\infty} p(n)\,t^n = (1 + t^1 + t^{1+1} + t^{1+1+1} + \cdots) \times (1 + t^2 + t^{2+2} + t^{2+2+2} + \cdots)$$
$$\times (1 + t^3 + t^{3+3} + t^{3+3+3} + \cdots) \times \cdots$$

(3.9.1)
$$= \prod_{n=1}^{+\infty} \frac{1}{1-t^n}.$$

See [15, Chap. 4] for much more about this.

The following proposition is the key to Formula (3.5.2) for an elliptic $K3$ surface $(S, L)$ as above. It yields Formula (3.5.2) after a series of manipulations of generating functions, involving Identity (3.9.1). I don't reproduce them here, and refer to [6, p. 383] for details. The rest of this subsection is dedicated to the proof of the proposition.

**(3.10) Proposition.** *The contribution of $\overline{M}_{\mathbf{a},\mathbf{b}}$ to $N_g^p$ is*

(3.10.1)
$$\left(\prod_{i=1}^{24} p(a_i)\right)\left(\prod_{j=1}^{g} b_j \sigma_1(b_j)\right).$$

**(3.11) Enumeration of elliptic multiple covers.** We first explain the factors $b_j \sigma_1(b_j)$ in (3.10.1); they are of a combinatorial nature.

Let $f : C \to S$ be a member of $\overline{M}_{\mathbf{a},\mathbf{b}}$. It is necessarily shaped as follows: the curve $C$ consists of (i) a smooth rational curve mapped isomorphically to the section $E \subseteq S$, which we will abusively call $E$ as well, (ii) $g$ smooth elliptic curves $G_1, \ldots, G_g$, pairwise disjoint and each attached at one point to $E$, and (iii) 24 trees of smooth rational curves $T_1, \ldots, T_{24}$, pairwise disjoint, each disjoint from the $G_j$'s and attached to $E$ at one point. For all $j = 1, \ldots, g$, $f$ maps $G_j$ to the elliptic fibre $F_j$ with degree $b_j$, and there is a marked point $x_j \in G_j$ mapped to $y_j$; for all $i = 1, \ldots, 24$, one has $f_* T_i = a_i R_i$.

For all $j$, we may fix the intersection point with $E$ as the origin of $G_j$ and $F_j$ respectively, which makes $f|_{G_j} : G_j \to F_j$ a degree $b_j$ homomorphism of elliptic curves. Such homomorphisms are in $1:1$ correspondence with index $b_j$ sublattices of the lattice defining $F_j$ as a complex torus, and the number of such sublattices is $\sigma_1(b_j)$ [39, § VII.5.2]. In addition, there are $b_j$ possibilities to choose the marked point $x_j$ in the preimage of $y_j$ in $G_j$.

Once the data of the homomorphisms $G_j \to F_j$ and the marked points $x_j \in G_f$ is fixed, the corresponding sub-moduli space of $\overline{M}_{\mathbf{a},\mathbf{b}}$ decomposes as a product

$$\prod_{i=1}^{24} \overline{M}_{a_i \mathbf{e}_i, 0},$$

where $\mathbf{e}_1 = (1, 0, \ldots, 0), \ldots, \mathbf{e}_{24} = (0, \ldots, 0, 1) \in \mathbf{Z}^{24}$. The moduli space $\overline{M}_{\mathbf{a},\mathbf{b}}$ thus consists of $\prod_{j=1}^{g} b_j \sigma_1(b_j)$ disjoint copies of the space $\prod_{i=1}^{24} \overline{M}_{a_i \mathbf{e}_i, 0}$, and the rest of the proof of Proposition (3.10) consists in showing that each of those contributes by $\prod_{i=1}^{24} p(a_i)$ to $N_g^p$.

**(3.12) Description of $\overline{M}_{a_i \mathbf{e}_i, 0}$.** Let us start by defining a model stable map $h_{a,R} : \mathbf{\Sigma}_a \to S$ for all positive integers $a$ and 1-nodal rational curves $R \in \{R_1, \ldots, R_{24}\}$. The curve $\mathbf{\Sigma}_a$ is a tree of $2a + 2$ smooth rational curves $\Sigma_E, \Sigma_{-a}, \ldots, \Sigma_0, \ldots \Sigma_{+a}$ as depicted on the figure below.

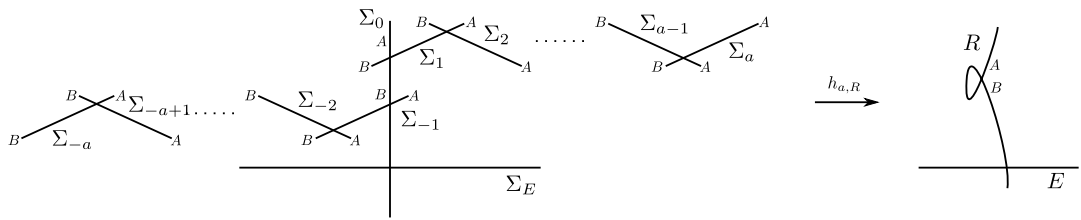

The map $h_{a,R}$ is chosen so that it restricts to an isomorphism $\Sigma_E \cong E$ (hence from now on we will denote $\Sigma_E$ by $E$) and to $2a+1$ copies $\Sigma_i \to R$ of the normalization of the 1-nodal rational curve $R$, in such a way that it is everywhere a local isomorphism between $\boldsymbol{\Sigma}_a$ and $E \cup R$. Concretely, the latter requirement is that locally at every node $\Sigma_i \cap \Sigma_{i+1}$, the map $h_{a,R}$ should send $\Sigma_i$ to one of the two local branches of $R$ at its node and $\Sigma_{i+1}$ to the other.

There are basically two possible (indifferent) choices. We indicate one of them on the above figure by decorating each local branch at a node of $\Sigma_{-a} \cup \ldots \cup \Sigma_{+a}$ with a letter $A$ or $B$, with the convention that $A$ and $B$ label the two local branches of $R$ at its node.

**(3.13) Lemma.** *For every $(f : C \to S) \in \overline{M}_{a\mathbf{e}_i,0}$, there is a unique lift of $f$ to a stable map $\tilde{f} : C \to \boldsymbol{\Sigma}_a$, meaning that $f = h_{a,R_i} \circ \tilde{f}$.*

*Proof.* The curve $C$ is necessarily made of a smooth rational component mapped isomorphically to $E$, which we denote by $E$ as well, and a tree of smooth rational curves $T$ attached at one point to $E$, as in (3.11). For each irreducible component $C_s$ of $T$, one has $f_*C_s = k_s R_i$ for some non-negative integer $k_s$. We determine the lift $\tilde{f}$ by exploring the dual graph of $T$ along all possible paths from its root to one of its leaves, as follows.

There is a unique irreducible component $C_0$ of $T$ intersecting $E$; we call the corresponding vertex of the dual graph of $T$ the *root* of the dual graph. The *leaves* are those vertices corresponding to irreducible components of $T$ intersecting only one other irreducible component. Now the lift $\tilde{f}$, should it exist, necessarily maps $C_0$ to $\Sigma_0$, and there is a unique suitable map $C_0 \to \Sigma_0$ (possibly contracting $C_0$ to the point $\Sigma_0 \cap E$) by the universal property of normalization.

Suppose a putative lift $\tilde{f}$ is determined on an irreducible component $C_s$ of $T$, and consider an arbitrary component $C_{s+1}$ of $T$ intersecting $C_s$ at one point $z_{s+1}$. We claim that the behaviour of $\tilde{f}$ on $C_{s+1}$ is uniquely determined by the already constructed piece $\tilde{f}|_{C_s}$. If $C_{s+1}$ is contracted by $f$, this is clear; otherwise it is enough, by the universal property of normalization, to determine which component of $\boldsymbol{\Sigma}_a$ the lift $\tilde{f}$ should map $C_{s+1}$ to. If $\tilde{f}(z_{s+1})$ is a smooth point of $\boldsymbol{\Sigma}_a$, then $\tilde{f}$ has to map $C_{s+1}$ to the same component it maps $C_s$ to; if not, $\tilde{f}(z_{s+1})$ is a node $\Sigma_{t_s} \cap \Sigma_{t_s+\varepsilon}$, $\varepsilon \in \{\pm 1\}$, and $\tilde{f}$ has to map $C_{s+1}$ to $\Sigma_{t_s}$ or $\Sigma_{t_s+\varepsilon}$, depending on which of the local branches $A$ and $B$ of $R_i$ at its node the local branch of $C_{s+1}$ at $z_{s+1}$ is mapped to by $f$.

This discussion shows that one may algorithmically construct an $\tilde{f}$ such that $h_{a,R_i} \circ \tilde{f} = f$, and that there is a unique such lift. (Note that the chain of rational curves $\Sigma_{-a} \cup \ldots \cup \Sigma_a$ is long enough for the construction to go through without any trouble: since $f_*C = E + aR_i$, if at some point during the algorithm we hit an irreducible component $C_s$ of $C$ that has to be mapped to $\Sigma_{-a+1}$ or $\Sigma_{a-1}$, then the push-forward by $f$ of the sum of all components already visited by the algorithm fills out the class $aR_i$, hence all components of $C$ not yet touched by the algorithm are contracted by $f$, and we don't have to go beyond $\Sigma_{-a+1}$ or $\Sigma_{a-1}$ in $\boldsymbol{\Sigma}_a$). □

Using Lemma (3.13), one can associate to every stable map $(f : C \to S) \in \overline{M}_{a\mathbf{e}_i,0}$ a combinatorial datum called an *admissible sequence* of weight $a$: this is a sequence of $2a+1$ non-negative integers $\mathbf{k} = (k_l)_{-a \leqslant l \leqslant a}$ with the property that there exist $m, n \geqslant 0$ such that $k_l = 0$ for all $l < -m$ and $l > n$, and $k_l > 0$ for all $l = -m, \ldots, n$. The association goes as follows: let $h_{a,R_i} \circ \tilde{f}$ be the factorization of $f$, and write

$$\tilde{f}_*C = E + \sum_{s=-a}^{a} k_s\Sigma_s.$$

It follows from the construction of $\tilde{f}$ in the proof of Lemma (3.13) that $(k_{-a},\ldots,k_a)$ is an admissible sequence of weight $a$. This gives the following.

**(3.14) Lemma.** *The moduli space $\overline{M}_{a\mathbf{e}_i,0}$ decomposes as the disjoint union*

$$\overline{M}_{a\mathbf{e}_i,0} = \coprod_{\mathbf{k}} \overline{M}_{\mathbf{k}},$$

*where $\mathbf{k}$ ranges over all weight $a$ admissible sequences, and $\overline{M}_{\mathbf{k}}$ is the sub-moduli space of $\overline{M}_{a\mathbf{e}_i,0}$ parametrizing those $f$ with associated admissible sequence $\mathbf{k}$.*

**(3.15) Identification of the virtual class.** It follows from the above that each moduli space $\overline{M}_{\mathbf{a},\mathbf{b}}$ decomposes as a disjoint union of various products $\prod_1^{24} \overline{M}_{\mathbf{k}_i}$. A key point in the proof of Bryan and Leung is the explicit identification of the restriction to the product $\prod_1^{24} \overline{M}_{\mathbf{k}_i}$ of the (reduced) virtual class giving rise to the invariant $N_g^p$ [6, § 5.2]. This is by far the most demanding part of their article, and I will not attempt to give any idea of the proof.

The upshot is that (i) the virtual class on $\prod_1^{24} \overline{M}_{\mathbf{k}_i}$ is a product of virtual classes on the various factors, and (ii) the virtual class on the factor $\overline{M}_{\mathbf{k}_i}$ is computed by means of a *virtual tangent bundle* on the target curve $\boldsymbol{\Sigma}_{a_i}$. This virtual tangent bundle $\mathsf{T}$ is the vector bundle on $\boldsymbol{\Sigma}_{a_i}$ defined by the conditions that it is isomorphic to $h_{a_i,R_i}^* T_S$ on $\Sigma_{-a} \cup \ldots \cup \Sigma_a$ and to $T_E \oplus \mathcal{O}_E(-1)$ on $E$.

Note that $h_{a_i,R_i}^* T_S$ restricts to $T_E \oplus \mathcal{O}_E(-2)$ on $E$; the correction made to define $\mathsf{T}$ corresponds to the fact that we want to kill the obstruction space $H^1(C, N_f)$ as we know the actual obstruction space is trivial although $H^1(C, N_f)$ is not, see (3.3).

**(3.16) Planar model.** Let $a$ be a non-negative integer. Thanks to the result discussed in (3.15), it is possible to construct a model for $\boldsymbol{\Sigma}_a$ embedded in a familiar surface $P$, such that the actual deformation theory of $\boldsymbol{\Sigma}_a$ in $P$ is isomorphic to the reduced virtual theory leading to the invariants $N_g^p$. This will eventually let us compute the contribution of the $\overline{M}_{\mathbf{k}_i}$'s (hence also of the $\prod_1^{24} \overline{M}_{\mathbf{k}_i}$'s) to the invariant $N_g^p$.

Let $p, p_{-1}, p_1$ be three distinct points lying on a line in the projective plane $\mathbf{P}^2$. Let $P_1 \to \mathbf{P}^2$ be the blow-up at these three points, and call $E$ the exceptional divisor over $p$, $\Sigma_0$ the proper transform of the line through $p, p_{-1}, p_1$, and $\Sigma_{-1}, \Sigma_1$ the exceptional divisors over $p_{-1}, p_1$ respectively. Next, recursively for all $s = 1, \ldots, a$, we perform the blow-up $P_{s+1} \to P_s$ at two general points of $\Sigma_{-s}$ and $\Sigma_s$ respectively, and let $\Sigma_{-s-1}$ and $\Sigma_{s+1}$ be the two corresponding exceptional divisors; we call $E, \Sigma_{-s}, \ldots, \Sigma_s$ the proper transforms in $P_{s+1}$ of the curves with the same name in $P_s$.

The curve $\Sigma_{-a} \cup \ldots \cup \Sigma_0 \cup \ldots \cup \Sigma_a$ in $P_{a+1}$ (note that this excludes the last two exceptional curves $\Sigma_{-a-1}$ and $\Sigma_{a+1}$) is isomorphic as an abstract curve to $\boldsymbol{\Sigma}_a$. Moreover, the tangent bundle of $P_{a+1}$ restricts to $\mathcal{O}_{\mathbf{P}^1}(2) \oplus \mathcal{O}_{\mathbf{P}^1}(-1)$ on $E$, and to $\mathcal{O}_{\mathbf{P}^1}(2) \oplus \mathcal{O}_{\mathbf{P}^1}(-2)$ on $\Sigma_{-a}, \ldots, \Sigma_a$; it is therefore isomorphic to the virtual tangent bundle $\mathsf{T}$ introduced in (3.15). Then, Bryan–Leung prove the following.

**(3.17) Lemma.** *For all admissible sequences $\mathbf{k} = (k_{-a}, \ldots, k_0, \ldots, k_a)$, the "local" contribution $\int_{[\overline{M}_{\mathbf{k}}]^{\mathrm{vir}}} 1$ of $\overline{M}_{\mathbf{k}}$ to the invariant $N_g^p$ equals the ordinary genus $0$ Gromov–Witten integral $\int_{[\overline{M}_0(P_{a+1},\beta)]^{\mathrm{vir}}} 1$ for the class $\beta = E + \sum_{-a}^{a} k_s \Sigma_s$.*

Thomas Dedieu                                                                                                    21

The lemma follows from (3.15) and the isomorphism between the restriction of $T_{P_{a+1}}$ and $\mathsf{T}$ to $\Sigma_a$, provided the two moduli spaces $\overline{M}_{\mathbf{k}}$ and the ordinary $\overline{M}_0(P_{a+1}, E + \sum_{-a}^a k_s \Sigma_s)$ are isomorphic as sets. Bryan and Leung are able to prove this by elementary arguments using a more elaborate set-up: they start with a linear $\mathbf{C}^*$ action (not $(\mathbf{C}^*)^2$) on $\mathbf{P}^2$ leaving the line $\langle p, p_{-1}, p_1 \rangle$ and a point $q$ fixed, and this $\mathbf{C}^*$ action survives in $P_{a+1}$. See [6, Lemma 5.7] for the proof.

Finally, by deformation invariance of Gromov–Witten invariants, we may transport the computation on the projective plane blown-up at $2a+3$ general points. This gives the following.

**(3.18) Lemma.** *Let $\tilde{\mathbf{P}}^2$ be the blow-up of $\mathbf{P}^2$ at a general set of $2a+3$ points, with exceptional divisors $E, E_{-a-1}, \ldots, E_{-1}, E_1, \ldots, E_{a+1}$ (all $(-1)$-curves of course). We call $H$ the pull-back of the line class. Then for all admissible sequences $\mathbf{k} = (k_{-a}, \ldots, k_0, \ldots, k_a)$, the "local" contribution $\int_{[\overline{M}_{\mathbf{k}}]^{\mathrm{vir}}} 1$ of $\overline{M}_{\mathbf{k}}$ to the invariant $N_g^p$ equals the ordinary genus $0$ Gromov–Witten integral $\int_{[\overline{M}_0(\tilde{\mathbf{P}}^2, \beta_{\mathbf{k}})]^{\mathrm{vir}}} 1$ for the class*

$$\beta_{\mathbf{k}} = E + \sum_{s=1}^a k_{-s}(E_{-s} - E_{-s-1}) + k_0(H - E - E_1 - E_{-1}) + \sum_{s=1}^a k_s(E_s - E_{s+1})$$

$$= k_0 H - (k_0 - 1)E - \sum_{s=1}^a (k_{s-1} - k_s)E_s - k_a E_{a+1} - \sum_{s=1}^a (k_{-s+1} - k_{-s})E_{-s} - k_{-a}E_{-a-1}.$$

Note that

$$(k_0 - 1) + \sum_{s=1}^a (k_{s-1} - k_s) + k_a + \sum_{s=1}^a (k_{-s+1} - k_{-s}) + k_{-a} = 3k_0 - 1,$$

so the virtual class $[\overline{M}_0(\tilde{\mathbf{P}}^2, \beta_{\mathbf{k}})]^{\mathrm{vir}}$ has dimension $0$, see (3.19) below.

**(3.19) The computation on the blown-up plane.** The Gromov–Witten invariants gotten in Lemma (3.18) are computable in practice thanks to the analysis of genus $0$ Gromov–Witten invariants of blow-ups of $\mathbf{P}^2$ carried out by Göttsche and Pandharipande [19].

Let $n$ be a non-negative integer, $d, \alpha_1, \ldots, \alpha_n$ integers. We call $N(d; \alpha_1, \ldots, \alpha_n)$ the genus $0$ Gromov–Witten invariant for $\tilde{\mathbf{P}}^2$ and the class $dH - \sum_i \alpha_i E_i$ (mind the minus sign, introduced for obvious geometric reasons), where $\tilde{\mathbf{P}}^2$ is the projective plane blown-up at a general set of $n$ points, $H$ the pull-back of the line class, and $E_1, \ldots, E_n$ the exceptional $(-1)$-curves. The corresponding moduli space of stable maps has virtual dimension $3d - 1 - \sum_i \alpha_i$; if this is positive, then we impose in addition the appropriate number of point constraints, and if this is negative, then we set the invariant to $0$.

These invariants enjoy the following properties (see [19]):
(i) $N(1) = 1$;
(ii) $N(d; \alpha_1, \ldots, \alpha_{n-1}, 1) = N(d; \alpha_1, \ldots, \alpha_{n-1}, 0) = N(d; \alpha_1, \ldots, \alpha_{n-1})$;
(iii) $N(d; \alpha_1, \ldots, \alpha_n) = N(d; \alpha_{\sigma(1)}, \ldots, \alpha_{\sigma(n)})$ for any permutation $\sigma \in \mathfrak{S}_n$;
(iv) $N(d; \alpha_1, \ldots, \alpha_n) = 0$ if there is an index $i$ for which $\alpha_i < 0$, unless $dH - \sum_i \alpha_i E_i = E_{i_0}$ for some $i_0$ in which case the invariant is $1$;
(v) The invariant $N(d; \alpha_1, \alpha_2, \alpha_3, \alpha_4, \ldots, \alpha_n)$ is invariant under the isomorphism given by the quadratic Cremona transformation corresponding to the linear system $|2H - E_1 - E_2 - E_3|$, i.e.,

$$N(d; \alpha_1, \alpha_2, \alpha_3, \alpha_4, \ldots, \alpha_n) =$$
$$N(2d - \alpha_1 - \alpha_2 - \alpha_3, d - \alpha_2 - \alpha_3, d - \alpha_1 - \alpha_3, d - \alpha_1 - \alpha_2, \alpha_4, \ldots, \alpha_n).$$

We need the following definition to state the result. An admissible sequence $(k_{-a}, \ldots, k_0, \ldots, k_a)$ is 1-*pyramidal*[5] if

$$k_s - 1 \leqslant k_{s+1} \leqslant k_s \quad \text{and} \quad k_{-s} - 1 \leqslant k_{-s-1} \leqslant k_{-s}$$

for all $s = 0, \ldots, a-1$.

**(3.20) Lemma.** *Let* $\mathbf{k} = (k_{-a}, \ldots, k_0, \ldots, k_a)$ *be an admissible sequence of weight* $a$. *Then the genus* 0 *Gromov–Witten invariant*

$$N(\mathbf{k}) := \int_{[\overline{M}_0(\tilde{\mathbf{P}}^2, \beta_{\mathbf{k}})]^{\mathrm{vir}}} 1$$
$$= N(k_0;\, k_0 - 1,\, k_{-a},\, k_{-a+1} - k_{-a},\, \ldots,\, k_0 - k_{-1},\, k_0 - k_1,\, \ldots,\, k_{a-1} - k_a,\, k_a)$$

*equals* 1 *if* $\mathbf{k}$ *is* 1-*pyramidal, and* 0 *otherwise.*

*Proof.* We first show that $N(\mathbf{k}) = 0$ if $\mathbf{k}$ is not 1-pyramidal. Since $k_0 > 0$ by definition of an admissible sequence, it follows from Property (iv) above that $N(\mathbf{k}) \neq 0$ implies

(3.20.1) $$k_{-s} \leqslant k_{-s+1} \quad \text{and} \quad k_{s-1} \geqslant k_s$$

for all $s = 1, \ldots, a$. Next, we apply for all $s = 2, \ldots, a-1$ the Cremona transformation defined by $|2H - E - E_1 - E_{s+1}|$ and get

$$N(\mathbf{k}) = N(1 + k_1 + k_{s+1} - k_s;\, \underbrace{k_1 - k_s + k_{s+1}}_{E},\, \ldots,\, \underbrace{1 - k_s + k_{s+1}}_{E_1},\, \ldots,\, \underbrace{1 - k_0 + k_1}_{E_{s+1}},\, \ldots)$$

by Property (v) above. If $N(\mathbf{k}) \neq 0$, we have $k_s \leqslant k_1$ by (3.20.1), hence $1 + k_1 + k_{s+1} - k_s > 0$, and Property (iv) then implies that $k_1 \geqslant k_0 - 1$ and $k_{s+1} \geqslant k_s - 1$. An analogous move shows that $k_{-s-1} \geqslant k_{-s} - 1$ for all $s = 0, \ldots, a-1$ if $N(\mathbf{k}) \neq 0$. Thus, the non-vanishing of $N(\mathbf{k})$ implies that $\mathbf{k}$ is 1-pyramidal.

Conversely, assume that $\mathbf{k}$ is 1-pyramidal of weight $a$. Then $k_a = k_{-a} = 0$ (otherwise the weight exceeds $a$; we have somehow already made this observation in the course of the proof of Lemma (3.13)), and all coefficients $k_s - k_{s-1}$ and $k_{-s} - k_{-s+1}$ for $s = 1, \ldots, a$, equal 0 or 1. It thus follows from Property (ii) that $N(\mathbf{k}) = N(k_0; k_0 - 1)$, which is readily seen to equal 1: the moduli space of genus 0 stable maps $\overline{M}_0(\tilde{\mathbf{P}}^2, k_0 H - (k_0 - 1)E)$ has only one enumeratively meaningful irreducible component, isomorphic to the family of degree $k_0$ plane curves with multiplicity $k_0 - 1$ at a fixed point $x_E \in \mathbf{P}^2$, and this is a linear system. $\square$

The proof of Theorem (3.5) will be complete once we show the following combinatorial result.

**(3.21) Lemma.** *The number of* 1-*pyramidal admissible sequences* $(k_{-a}, \ldots, k_0, \ldots, k_a)$ *of weight* $a$ *equals the partition number* $p(a)$ *(cf. (3.9)).*

Indeed, together with Lemmata (3.14), (3.18) and (3.20), this shows that the "local" contribution of $\overline{M}_{a_i \mathbf{e}_i, 0}$ equals $p(a_i)$, hence the contribution of each copy of the product $\prod_1^{24} \overline{M}_{a_i \mathbf{e}_i, 0}$ equals $\prod_1^{24} p(a_i)$; this proves Proposition (3.10), thanks to the enumeration of elliptic multiple covers performed in (3.11). Finally, as already mentioned, Proposition (3.10) implies Theorem (3.5).

---

[5]this is called 1-*admissible* in [6].

*Proof of Lemma (3.21).* Partitions of an integer $a$ are in bijective correspondence with Young diagrams of size $a$ [15, Chap. 4]. We shall exhibit a bijective correspondence between Young diagrams of size $a$ and 1-pyramidal admissible sequences $(k_{-a}, \ldots, k_0, \ldots, k_a)$ of weight $a$. We see Young diagrams as embedded in the upper-right quadrant of a Cartesian plane, leaning on both the $x$ and $y$ axes, and with blocks squares of size 1. Given such a Young diagram, we let $k_s$ be the number of blocks on the line $y - x = s$ for all $s = -a, \ldots, 0, \ldots, a$: see the figure below for an example of this procedure.

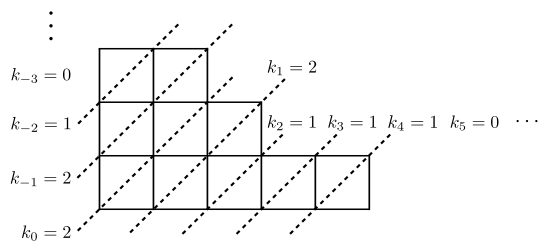

We leave it to the reader to check that $(k_{-a}, \ldots, k_a)$ is indeed admissible, 1-pyramidal, of weight $a$, and that this correspondence establishes the announced bijection. □

## 4 – BPS state counts

In this Section, I discuss why and how curve counting in non-primitive classes imply the use of multiple covers formulæ. This features the generalization of the Yau–Zaslow formula of Theorem (2.1) to non-primitive classes.

### 4.1 – Rational curves on the quintic threefold

To describe the picture in its simplest form, let us first discuss a question that lies slightly outside the scope of these notes, namely that of counting rational curves on a general quintic hypersurface $V$ of $\mathbf{P}^4$.

**(4.1)** The Clemens conjecture (see (4.5) below) predicts that there are finitely many such curves of any given degree $d$, in keeping with the virtual dimension

$$(4.1.1) \qquad \operatorname{vdim} \overline{M}_g(V, \beta) = (\dim V - 3)(g - 1) - K_V \cdot \beta$$

being 0 for any homology class on a Calabi–Yau threefold $V$. This suggests that the number

$$N_d^V := \int_{[\overline{M}_0(V, d\ell)]^{\mathrm{vir}}} 1$$

(where $\ell$ denotes the homology class of a line) may give the actual number of rational curves of degree $d$ in $V$. This would be particularly appealing, since the numbers $N_d^V$ may in theory be rigorously computed using the predictions of mirror symmetry; this has been done in practice for small values of $d$, see [10] for a thorough discussion.

A first objection against such an ideal correspondence is that there may be rational curves with non-trivial infinitesimal deformations, but the Clemens conjecture also predicts that this does not happen.

**(4.2) Multiple cover formula.** An actual issue preventing $N_d^V$ from directly reflecting the number of degree $d$ rational curves on $V$ is the existence of components of $\overline{M}_0(V, d\ell)$ having dimension larger than expected (i.e., positive dimension). Suppose there is a smooth degree $d_0 = d/k$ rational curve $\Gamma \subseteq V$ for some integer $k > 1$, and let $f_0 : \mathbf{P}^1 \to V$ be the stable map induced by the normalization of $\Gamma$. Then, for all degree $k$ covers $\phi_k : \mathbf{P}^1 \to \mathbf{P}^1$, the composed map $f_0 \circ \phi_k$ is a member of $\overline{M}_0(V, d\ell)$, and the locus in $\overline{M}_0(V, d\ell)$ corresponding to such composed maps over $\Gamma$ is an irreducible variety $M_k(\Gamma)$ of dimension $2k - 2$.

The Aspinwall–Morrison formula asserts that the corresponding irreducible component $\overline{M}_k(\Gamma)$ of $\overline{M}_0(V, d\ell)$ contributes by $\frac{1}{k^3}$ to the integral $N_d^V$. This has been mathematically proved by Kontsevich and Manin, and Voisin, see [43, § 5.6] and [10, Thm. 7.4.4].

To make this plausible it is convenient to consider, instead of the integrals $N_d^V$, their close cousins
$$\langle I_{0,3,d\ell}^V \rangle(\omega_1, \omega_2, \omega_3) := \int_{[\overline{M}_{0,3}(V,d\ell)]^{\mathrm{vir}}} \mathrm{ev}_1^*(\omega_1) \wedge \mathrm{ev}_2^*(\omega_2) \wedge \mathrm{ev}_3^*(\omega_3),$$
where $\overline{M}_{0,3}(V, d\ell)$ is the space of genus 0 stable maps with 3 marked points, and the $\omega_i$ are Kähler forms on $V$. The moduli space $\overline{M}_{0,3}(V, d\ell)$ has the advantage over $\overline{M}_0(V, d\ell)$ that it identifies with the Hilbert scheme $\mathrm{Hom}(\mathbf{P}^1, V)$, locally at those stable maps having source $\mathbf{P}^1$. Now, it follows from the divisorial axiom of Gromov–Witten theory that
$$\langle I_{0,3,d\ell}^V \rangle(\omega_1, \omega_2, \omega_3) = \int_{d\ell} \omega_1 \times \int_{d\ell} \omega_2 \times \int_{d\ell} \omega_3 \times N_d^V.$$
On the other hand, each physical degree $d$ rational curve $C \subseteq V$ contributes via its normalization $f : \mathbf{P}^1 \to V$ by
$$\int_{\mathbf{P}^1} f^*\omega_1 \times \int_{\mathbf{P}^1} f^*\omega_2 \times \int_{\mathbf{P}^1} f^*\omega_3$$
to the integral $\langle I_{0,3,d\ell}^V \rangle(\omega_1, \omega_2, \omega_3)$, hence by
$$\frac{\left(\int_{\mathbf{P}^1} f^*\omega_1\right)\left(\int_{\mathbf{P}^1} f^*\omega_2\right)\left(\int_{\mathbf{P}^1} f^*\omega_3\right)}{\left(\int_{d\ell} \omega_1\right)\left(\int_{d\ell} \omega_2\right)\left(\int_{d\ell} \omega_3\right)}$$
to the invariant $N_d^V$. The contribution of $\Gamma$ to $N_{d_0}^V$ may be expressed in the same way.

Now, suppose we knew what the compactification $\overline{M}_k(\Gamma)$ of $M_k(\Gamma)$ in $\overline{M}_{0,3}(V, d\ell)$ looks like, and we had a vector bundle $E$ over it having fibre over $(f_0 \circ \phi_k) \in M_k(\Gamma)$ the obstruction space $H^1(\mathbf{P}^1, (f_0 \circ \phi_k)^* T_V)$. Then we could compute the contribution of $\overline{M}_k(\Gamma)$ to $\langle I_{0,3,d\ell}^V \rangle(\omega_1, \omega_2, \omega_3)$ by the excess formula

(4.2.1) $$\int_{\overline{M}_k(\Gamma)} c_{2k-2}(E) \wedge \mathrm{ev}_1^*(\omega_1) \wedge \mathrm{ev}_2^*(\omega_2) \wedge \mathrm{ev}_3^*(\omega_3).$$

The computation of [43, p. 115–116] shows that a heuristic model for $\overline{M}_k(\Gamma)$ leads to the following expression for this excess contribution:
$$\int_{\mathbf{P}^1} f_0^*\omega_1 \times \int_{\mathbf{P}^1} f_0^*\omega_2 \times \int_{\mathbf{P}^1} f_0^*\omega_3.$$
The upshot is that the direct contribution of $\Gamma$ to $N_{d_0}^V$ and its contribution via $\overline{M}_k(\Gamma)$ to $N_d^V$ (in this heuristic model) are, respectively,
$$\frac{\left(\int_{\mathbf{P}^1} f_0^*\omega_1\right)\left(\int_{\mathbf{P}^1} f_0^*\omega_2\right)\left(\int_{\mathbf{P}^1} f_0^*\omega_3\right)}{\left(\int_{d_0\ell} \omega_1\right)\left(\int_{d_0\ell} \omega_2\right)\left(\int_{d_0\ell} \omega_3\right)} \quad \text{and} \quad \frac{\left(\int_{\mathbf{P}^1} f_0^*\omega_1\right)\left(\int_{\mathbf{P}^1} f_0^*\omega_2\right)\left(\int_{\mathbf{P}^1} f_0^*\omega_3\right)}{\left(\int_{d\ell} \omega_1\right)\left(\int_{d\ell} \omega_2\right)\left(\int_{d\ell} \omega_3\right)}.$$
The latter is $\frac{1}{k^3}$ times the former, since $\int_{d\ell} \omega_i = k \int_{d_0\ell} \omega_i$ for all $i = 1, 2, 3$.

Thomas Dedieu                                                                                          25

**(4.3) Instanton numbers.** The idea is to define new invariants $n_d^V$ from the Gromov–Witten integrals $N_d^V$ that better reflect the enumerative geometry of the Calabi–Yau threefold $V$, taking into account the multiple cover phenomenon described in paragraph (4.2) above. The relations that these numbers should satisfy are

(4.3.1) $$N_d^V = \sum_{k|d} \tfrac{1}{k^3} n_{d/k}^V$$

for all positive integers $d$, where the sum runs over all positive integral divisors of $d$. They form an invertible triangular set of relations, and each number $n_d^V$ is therefore uniquely determined by the invariants $N_{d'}^V$ for all positive divisors $d'$ of $d$. These numbers $n_d^V$ are traditionally called *instanton numbers*. This name obviously bears some physical meaning, but let me not try to discuss it.

**(4.4) Reducible covers of singular curves.** There is yet another phenomenon, first observed by Pandharipande, that prevents the instanton numbers $n_d^V$ defined in (4.3) above from being the actual numbers of integral degree $d$ rational curves on $V$. It comes from the existence of singular integral rational curves.

To describe the simplest instance of this phenomenon, assume there is a degree $d$ integral rational curve $C \subseteq V$ with an ordinary double point $x$. Let $f : \mathbf{P}^1 \to V$ be the normalization of $C$, and denote by $x_1$ and $x_2$ the two points in the preimage of $x$ in the normalization, Now, consider the map
$$f_x : \mathbf{P}^1 \cup_{x_1,x_2} \mathbf{P}^1 \to V,$$
obtained by glueing tranversaly two copies of the normalization of $C$, in such a way that $x_1$ in the first copy is identified with $x_2$ in the second copy. It is a stable map of genus 0 realizing the homology class $2[C]$ on $V$, hence it contributes to the Gromov–Witten invariant $N_{2d}^V$. This contribution is by 1 if $C$ is infinitesimally rigid, cf. [10, § 9.2.3]. It follows that a $\delta$-nodal rational curve of degree $d$ (i.e., a curve with exactly $\delta$ ordinary double points as singularities) contributes in the above described fashion by $\delta$ to the Gromov–Witten integral $N_{2d}^V$.

The Clemens conjecture below predicts that the complications don't go beyond this in the case of a general quintic threefold. Note in particular that part (iii) of the conjecture implies that the numbers $N_d^V$ (or $n_d^V$) count *irreducible* physical rational curves $C \subseteq V$, since by definition the source of a stable map is connected.

**(4.5) Conjecture** (Clemens). *Let $V \subseteq \mathbf{P}^4$ be a general quintic hypersurface.*
*(i) For each integer $d \geqslant 1$, there are only finitely many irreducible rational curves $C \subseteq V$ of degree $d$.*
*(ii) For every integral rational curve $C \subseteq V$ with normalisation $f : \mathbf{P}^1 \to V$, the normal bundle $N_f$ of the map $f$ is isomorphic to $\mathcal{O}_{\mathbf{P}^1}(-1) \oplus \mathcal{O}_{\mathbf{P}^1}(-1)$ (i.e., $C$ is infinitesimally rigid).*
*(iii) All the integral rational curves on $V$ (of any degree) are pairwise disjoint.*

It is completely proved in degree $\leqslant 11$ (cf. [9, 11] for the latest steps), and part (i) is known in degree 12 [1].

We end this prologue about quintic threefolds by an explicit example displaying all these phenomena together.

**(4.6) Rational curves of degree 10 on a quintic threefold.** (cf. [10, § 9.2.3] for a thorough analysis). First note that, by definition of the instanton numbers, the Gromov–Witten integral $N_{10}^V$ decomposes as
$$N_{10} = \frac{1}{10^3} n_1 + \frac{1}{5^3} n_2 + \frac{1}{2^3} n_5 + n_{10}$$

(we have now dropped the superscript $V$ to lighten notations), and one has

$$n_1 = 2,875$$
$$n_2 = 609,250$$
$$n_5 = 229,305,888,887,625$$
$$n_{10} = 704,288,164,978,454,686,113,488,249,750,$$

cf. [7]. While $n_1$ and $n_2$ are simply the numbers of lines and conics on a general quintic threefold, $n_5$ counts two kinds of rational curves of degree 5. Indeed, the planes in $\mathbf{P}^4$ are parametrized by a 6-dimensional Grassmannian, and for a general quintic $V \subseteq \mathbf{P}^4$, finitely many of them are 6-tangent to $V$; the corresponding plane sections of $V$ are 6-nodal plane quintic curves, and in particular they are rational curves. Vainsencher [41] has been able to compute their number

$$n'_5 = 17,601,000,$$

see [\*\*, Paragraph \*\*]. Each such curve contributes to $n_5$ (or $N_5$) by 1 [10, Lemma 9.2.4], and

$$n''_5 := n_5 - n'_5$$

is indeed the number of *smooth* rational curves of degree 5 on $V$.

Whereas $n_5$ still is the number of rational curves of degree 5 on $V$, this is no longer true for $n_{10}$ as the discussion in (4.4) indicates. The actual number of degree 10 integral rational curves on $V$ is

$$n^\circ_{10} := n_{10} - 6n'_5$$

[10, Thm. 9.2.6][6].

## 4.2 – Degree $8$ rational curves on a sextic double plane

In this subsection we give the enumerative interpretation due to Gathmann [17] of the reduced Gromov–Witten invariant

$$N^2_{0,2} := \int_{[\overline{M}_0(S,2L)]^{\mathrm{red}}} 1$$

of a general primitively polarized $K3$ surface $(S, L)$ of genus 2 (i.e., $S$ is a double covering of the plane $\pi : S \to \mathbf{P}^2$ branched over a general sextic curve $B$, and $L$ is the pull-back of the line class). Note that the polarization $2L$ has genus 5, as $(2L)^2 = 8$.

**(4.7)** The analysis carried out in Section 4.1 indicates that the integral $N^2_{0,2}$ is a sum of contributions corresponding to the following types of curves.
*(i)* 5-nodal integral curves in $|2L|$; these are the preimages of the conics in $\mathbf{P}^2$ tangent to the branch sextic $B$ at 5 distinct points. There are $70,956$ of those, as Gathmann was able to compute using his theory of relative Gromov–Witten invariants for hypersurfaces.
*(ii)* Reducible rational curves made of two distinct rational curves in $|L|$: let $C_1, C_2 \in |L|$ be two rational curves (each of these is the pull-back of a line bitangent to $B$); they intersect transversely in two points $x, x'$ where both of them are smooth, if $(S, L)$ is general; then, the partial normalization of $C_1 + C_2$ at the four nodes of $C_1$ and $C_2$ and one node among $x$ and $x'$ gives a stable map with source a chain of two smooth rational curves, mapping birationally onto the physical curve $C_1 + C_2 \in |2L|$.

These two stable maps contribute each by 1 to $N^2_{0,2}$, so in the upshot each of the $\binom{324}{2}$ pairs of distinct rational curves in $|L|$ contributes by 2 in this way, thus giving a total contribution of $104,652$ to $N^2_{0,2}$ (recall that there are 324 bitangent lines to $B$, as can be classicaly computed, or extracted from Theorem (2.1)).

---

[6]there is a misprint there: $6\frac{1}{8}$ should be replaced by $6 + \frac{1}{8}$.

*(iii)* Reducible double coverings of one rational curve in $|L|$, as in (4.4). Since all rational curves in $|L|$ are 2-nodal, all 324 of them give 2 stable maps with reducible source, contributing by 1 each to $N^2_{0,2}$, for a total contribution of 648.

*(iv)* Irreducible double covers of one rational curve in $|L|$, as in (4.2). Gathmann [17, Lemma 4.1] proved that each of the 324 corresponding positive-dimensional irreducible components of $\overline{M}_0(S, 2L)$ gives a contribution by $\frac{1}{8}$ to $N^2_{0,2}$, similar to the prediction of the Aspinwall–Morrison formula, although we are in a significantly different situation (we are not on a Calabi–Yau threefold, and the rational curves we are considering double covers of are singular, for instance).

**(4.8) Remark.** An important difference with the case of the quintic threefold, although rather innocent-looking, is that in the present situation integral rational curves do intersect each other, contrary to the prediction of part (iii) of Clemens' Conjecture. In the above discussion (4.7), this amounts to case (ii) needing to be added with respect to the discussion for quintic threefolds, which is still manageable. When looking at linear systems $|mL|$ with $m \geqslant 3$ however, this soon gets much more complicated, see (4.15).

**(4.9)** From the enumeration of (4.7), one may deduce the value of $N^2_{0,2}$, which was not known before. But the striking observation of [17] is that the sum of the contributions (i)–(iii) above, which would be the instanton number $n^2_{0,2}$ in the language of (4.3), actually equals

$$N^5_{0,1} := N^5_0 = 176,256,$$

the number of degree 8 rational curves in a primitively polarized $K3$ surface of genus 5, computed by formula (2.1.1). This suggests the amazing possibility that the number of rational curves in a linear system $|L|$ on a $K3$ surface is independent of the divisibility of the line bundle $L$. We will come back to this in detail in Section 4.4 below.

## 4.3 – Elliptic curves in a 2-divisible class on a $K3$ surface

Here, I report on the computation by Lee and Leung [24] of the reduced Gromov–Witten invariant

$$N^p_{1,2} := \int_{[\overline{M}_{1,1}(S,2L)]^{\mathrm{red}}} \mathrm{ev}^*_1(\mathrm{pt})$$

of a general primitively polarized $K3$ surface $(S, L)$ of genus $p$, which "counts" genus 1 curves in $|2L|$ passing through 1 general point on $S$. Using a suitable version of topological recursion, they prove the following formula.

**(4.10) Theorem.** [24] *One has*

$$N^p_{1,2} = N^{4p-3}_{1,1} + 2N^p_{1,1}.$$

The numbers $N^q_{1,1} := N^q_1$ are those giving the number of elliptic curves through a general point in the primitive class of a $K3$ surface of genus $q$, as in formula (3.5). Note that $p' := 4p - 3$ is the genus of the polarization $2L$, i.e., $(2L)^2 = 2p' - 2$.

**(4.11)** Lee and Leung propose the following interpretation of their formula (4.10), to put it in tune with (4.9), and more generally with the results of Section 4.4 below.

Given a smooth elliptic curve $E$, there are $\sigma_1(2) = 1 + 2 = 3$ morphisms of elliptic curves $E' \to E$ of degree 2 (note that we require that the origin is respected), as we have already seen in (3.11). This implies that each of the $N^p_{1,1}$ elliptic curves $C$ in the primitive linear system $|L|$ passing through a general point $x_1 \in S$ gives via double covers of its normalization three genus 1 stable maps realizing the homology class $2[C]$, each contributing by 1 to the number $N^p_{1,2}$.

Lee and Leung deduce from this that $N^{4p-3}_1$ is the actual number of physical elliptic curves in $|2L|$, meaning that it counts each integral elliptic curve $C \in |2L|$ for 1 and each $2C$ where $C \in |L|$ is an integral elliptic curve for 1 as well. This is indeed a striking interpretation, although arguably debatable.

There is at any rate a phenomenon that prevents this interpretation to be anything more than philosophical, namely that reducible curves $C_0 + C_1$, where $C_0$ (resp. $C_1$) is a rational (resp. elliptic) integral curve in $|L|$, also contribute to the invariant $N^p_{1,2}$. The two curves $C_0$ and $C_1$ intersect (transversely, say) at $2p-2$ points $y_1, \ldots, y_{2p-2}$, and this gives $2p-2$ genus 1 stable maps
$$\mathbf{P}^1 \cup_{y_i} \bar{C}_1 \to S$$
realizing the class $2L$ and passing through the appropriate fixed point whenever $C_1$ does (the source is the transverse union of the normalization $\mathbf{P}^1$ and $\bar{C}_1$ of $C_0$ and $C_1$, attached at one point lying in the preimages of one of the points $y_i$).

## 4.4 – The Yau–Zaslow formula for non-primitive classes

We now come to the general statement confirming (4.9) and recently proved by Klemm, Maulik, Pandharipande, and Scheidegger.

**(4.12)** Let $S$ be a $K3$ surface, and $L \in \mathrm{Pic}(S)$. It follows from deformation invariance of reduced Gromov–Witten integrals and the global Torelli theorem for $K3$ surfaces that the integral $\int_{[\overline{M}_0(S,L)]^{\mathrm{red}}} 1$ only depends on the self-intersection $L^2$ and the divisibility index of $L$ in $\mathrm{Pic}(S)$ (the latter is the largest integer $m$ for which there exists $L' \in \mathrm{Pic}(S)$ such that $L = mL'$). We may thus make the following definition.

For all integers $p > 0$ and $m \geqslant 1$, we let
$$N^p_{0,m} := \int_{[\overline{M}_0(S,mL)]^{\mathrm{red}}} 1$$
where $(S, L)$ is any primitively polarized $K3$ surface of genus $p$.

**(4.13) BPS states.** Similar to what has been done in (4.3), and following the insight of (4.9), we now define new invariants from the $N^p_{0,m}$ of (4.12) above by applying the corrections indicated by the Aspinwall–Morrison formula. Let us formulate this in terms of generating series as follows. Given a positive integer $p$, set

(4.13.1) $$F^p(v) := \sum_{m \geqslant 1} N^p_{0,m} v^m$$

as a formal power series in the variable $v$. Then the new set of invariants $n^p_{0,m}$ is uniquely determined by the rewriting of the generating series as

(4.13.2) $$F^p(v) := \sum_{m \geqslant 1} n^p_{0,m} \left( \sum_{d > 0} \frac{1}{d^3} v^{dm} \right)$$

(this is exactly the same modification as that of (4.3.1)).

Note that this is not mere makeshift reformulation. The invariants $n^p_{0,m}$ are believed to count objects named BPS states by the physicists, after Bogomol'nyi, Prasad, and Sommerfield, and have been introduced by Gopakumar and Vafa. The mathematical nature of these BPS states is however not clear yet. In particular, it should be possible to define the $n_{0,m}^p$ intrinsically, without relying on the $N_{0,m}^p$; the relation (4.13.2) would then tie together these two sets of independently defined invariants. See the enlightening survey [32, § $2\frac{1}{2}$] for more about this. There is moreover a physical meaning to the introduction of generating series, but I will not discuss it.

**(4.14) Theorem** ([23]). *The genus $0$ invariants $n_{0,m}^p$ do not depend upon the divisibility index, i.e., for all integers $m, p \geqslant 1$, one has*

(4.14.1) $$n_{0,m}^p = n_{0,1}^{m^2p-m^2+1} = N_0^{m^2p-m^2+1}$$

*(the integer $p' = m^2p - m^2 + 1$ is designed such that $(mL)^2 = 2p' - 2$ if $L^2 = 2p - 2$).*

Recall that $N_0^p$ was defined in section 2; the second equality in (4.14.1) is by definition of $n_{0,1}^p$ and $N_0^p$. This statement was part of the Yau–Zaslow conjecture [44]. Together with Theorem (2.1), which was also part of the Yau–Zaslow conjecture, it implies that all the $n_{0,m}^p$'s may be computed by means of formula (2.1.1). The set of relations (4.13.1) being triangular invertible, this also gives all genus 0 reduced Gromov–Witten invariants of $K3$ surfaces. Section 6 below contains an overview of the proof given by Klemm, Maulik, Pandharipande, and Scheidegger [23] of Theorem (4.14).

The truly remarkable feature of the invariants $n_{0,m}^p$ displayed by this statement is that the number of rational curves of prescribed degree in an algebraic $K3$ surface does not depend on the algebraic geometry of the surface.

**(4.15)** In spite of formula (4.13.2) taking into account the Aspinwall-Morrison multiple cover correction, the invariants $n_{0,m}^p$ do not in general count the actual number of rational curves in $|mL|$.

One reason for this is the existence of more non-reduced curves with rational support than those taken in consideration in the correction (4.13.2), namely curves with reducible support. For instance let $m = 3$ and consider two integral rational curves $C_1, C_2$. Then $2C_1 + C_2 \in |3L|$, and there are correspondingly finitely many positive-dimensional components of $\overline{M}_0(S, 3L)$, the general points of which correspond to stable maps

$$\mathbf{P}^1 \cup_x \mathbf{P}^1 \to S$$

with source a transverse union of two $\mathbf{P}^1$'s, consisting of a double cover of $C_1$ on the first component and the normalization of $C_2$ on the other. These certainly give an excess contribution to the invariant $N_{0,3}^p$, which is not taken into account in the definition of $n_{0,3}^p$.

Such problematic phenomena do not occur for $m \leqslant 2$, so that $n_{0,1}^p$ is directly enumerative as was already noted in Section 2, and $n_{0,2}^p$ counts reduced rational curves in the way described in Section 4.2. It would be very interesting to relate $n_{0,m}^p$ to the number of integral rational curves in $|mL|$ for $m \geqslant 3$.

There were, at least conjecturally, no such phenomena at work in the case of the quintic threefold discussed in subsection 4.1, as part (iii) of the Clemens conjecture (4.5) asserts that two integral rational curves in a general quintic threefold never intersect. In a surface, there is of course not enough space for two curves to avoid each other, so we inevitably have to deal with the aforementoined degenerate contributions.

The philosophy, as R. Pandharipande communicated to me, is that what the BPS numbers for $K3$ surfaces are virtually counting, are rational curves in some perturbation of the twistor family of the $K3$ surface (a threefold, cf. (3.2)). We shall consider in more detail the close interplay between counting invariants for $K3$ surfaces and Calabi–Yau threefolds in the next Section 5.

# 5 – Relations with threefold invariants

It is already visible in the very foundation of the theory of Gromov–Witten invariants for algebraic $K3$ surfaces developed by Bryan–Leung [6], see (3.2), that these invariants are fundamentally attached to a threefold (even though the approach by Maulik–Pandharipande [28], see (3.3), bypasses this). Another revealing evidence of the 3-dimensional nature of these invariants is the meaningful role played by the Aspinwall–Morrison in the Yau–Zaslow statement discussed in subsection 4.4 above, a tool specifically designed for Calabi–Yau threefolds. In this section we will try to describe this relation in a more conceptual way. It it wise to keep in mind the symplectic nature of Gromov–Witten invariants throughout.

## 5.1 – Two obstruction theories

**(5.1) A threefold degenerate contribution.** Let $V$ be a Calabi–Yau threefold. It follows from formula (4.1.1) that the virtual dimension of all spaces of stable maps is 0 regardless of the genus. This makes the following phenomenon happen.

Let $C_0 \subseteq V$ be a rational curve (smooth and infinitesimally rigid, say). Its normalization $f : \mathbf{P}^1 \to V$ contributes regularly by 1 to the integral $\int_{[\overline{M}_0(V,[C_0])]^{\mathrm{vir}}} 1$. But for any stable curve $C'$ of genus $g \geqslant 1$ we may obtain a genus $g$ stable map realizing the class $[C_0]$ by attaching $C'$ to the normalization of $C_0$ over a smooth point $x$, and letting

$$f_{C',x} : \mathbf{P}^1 \cup_x C_0 \to V$$

equal to $f$ along $\mathbf{P}^1$ and collapsing $C'$ to $x$. This produces a positive dimensional moduli space of genus $g$ stable maps which all have the same image $C_0 \subseteq V$; its contribution to the Gromov–Witten invariant $\int_{[\overline{M}_g(V,[C_0])]^{\mathrm{vir}}} 1$ must be computed via Hodge integrals over the moduli space of stable curves of genus $g$, see below. This has been studied by Faber and Pandharipande, see [32, § $1\frac{1}{2}$] and the references therein.

**(5.2) Curves in the twistor space of a $K3$.** Let $S$ be a $K3$ surface, together with an algebraic class $\beta \in H_2(S, \mathbf{Z})$. We consider its twistor space $T \to \mathbf{S}^2$ described in (3.2) above (we emphasize that this is a real 6-dimensional variety), and denote by $\iota : S \hookrightarrow T$ the canonical inclusion of $S$. Since curves in $T$ can only appear in the fibre $S$, we have the equality of moduli spaces of stable maps

$$\overline{M}_g(T, \iota_*\beta) = \overline{M}_g(S, \beta),$$

a priori only as sets but in fact as Deligne–Mumford stacks. They come however with two different obstructions theories, hence also with two different virtual classes. Gromov–Witten invariants on $T$ are related to those on $S$ (within the reduced theory for $K3$ surfaces, cf. (3.3)) by the formula

$$(5.2.1) \qquad \int_{[\overline{M}_g(T,\iota_*\beta)]^{\mathrm{vir}}} 1 \;=\; \int_{[\overline{M}_g(S,\beta)]^{\mathrm{red}}} (-1)^g \lambda_g,$$

where $\lambda_g$ stands for the top Chern class $c_g(\mathbf{E}_g)$ of the Hodge bundle $\mathbf{E}_g \to \overline{M}_g(S, \beta)$, whose fibre over the stable map $f : C \to S$ is $H^0(C, \omega_C)$.

**(5.3) Hodge integrals.** It follows from the invariance of reduced Gromov–Witten invariants under algebraic deformation and the global Torelli theorem for $K3$ surfaces that the right-hand side of (5.2.1) depends only on the self-intersection $\beta^2$ and the divisibility index of $\beta$ as an algebraic class. We may thus formulate the following definition.

For integers $g \geqslant 0$ and $p, m \geqslant 1$, let

$$(5.3.1) \qquad R_{g,m}^p := \int_{[\overline{M}_g(S, mL)]^{\mathrm{red}}} (-1)^g \lambda_g$$

where $(S, L)$ is any primitively polarized $K3$ surface of genus $p$, and $\lambda_g$ is the top Chern class of the Hodge bundle as in (5.2) above.

This extends the definition of the invariants $N_{0,m}^p$ in (4.12) above, in the sense that $N_{0,m}^p = R_{0,m}^p$ (note however that the invariants $N_{1,2}^p$ used in subsection 4.3 do not coincide with the $R_{1,2}^p$). For $g > 0$, these invariants are certainly not counting curves on $S$; rather, formula (5.2.1) tells us that they virtually give the excess contribution of $S$ to the vertical Gromov–Witten theory of any $K3$-fibred threefold in which it appears as a fibre. This philosophy is put into concrete form by Theorem (6.9) below.

## 5.2 – The Katz–Klemm–Vafa formula

This is an extension of the Yau–Zaslow conjecture discussed in Section 4.4 above to the invariants $R_{g,m}^p$. It has been proved by Pandharipande and Thomas [36], see also [35].

**(5.4) BPS invariants.** It is admittedly better to organize the invariants $R_{g,m}^p$ in BPS form as in (4.13). We have now a clear justification for this, as we have seen in Section 5.1 above that these invariants really count objects on threefolds.

We first let

$$F^p(u, v) := \sum_{g \geqslant 0} \sum_{m > 0} R_{g,m}^p u^{2g-2} v^m$$

as formal power series in the two variables $u, v$ for all positive integers $p$. One then defines new invariants $r_{g,m}^p$ for all integers $p, m > 0$, $g \geqslant 0$, by setting

$$\begin{aligned}
F^p(u, v) &= \sum_{m>0} \sum_{g \geqslant 0} r_{g,m}^p u^{2g-2} \left( \sum_{d>0} \frac{1}{d} \left( \frac{\sin d\frac{u}{2}}{\frac{u}{2}} \right)^{2g-2} v^{dm} \right) \\
&= \sum_{m>0} \Bigg( r_{0,m}^p u^{-2} \sum_{d>0} \Big( \frac{1}{d^3} + \frac{1}{12d} u^2 + \frac{d}{240} u^4 + \frac{d^3}{6048} u^6 + \frac{d^5}{172800} u^8 + \cdots \Big) v^{dm} \\
&\quad + r_{1,m}^p \sum_{d>0} \frac{1}{d} v^{dm} \\
&\quad + r_{2,m}^p u^2 \sum_{d>0} \Big( d - \frac{d^3}{12} u^2 + \frac{d^5}{360} u^4 - \frac{d^7}{20160} u^6 + \frac{d^9}{1814400} u^8 + \cdots \Big) v^{dm} \\
&\quad + r_{3,m}^p u^4 \sum_{d>0} \Big( d^3 - \frac{d^5}{6} u^2 + \frac{d^7}{80} u^4 - \frac{17 d^9}{30240} u^6 + \frac{31 d^{11}}{1814400} u^8 + \cdots \Big) v^{dm} \\
&\quad + \cdots \Bigg).
\end{aligned}$$

The modifications for genus $g > 0$ objects did not appear earlier in this text. Note that every object counted by $r_{g,m}^p$ contributes to the invariants $R_{g',m}^p$ for all $g' \geqslant g$ (except when $g = 1$), with alternated sign if $g \geqslant 2$. This is in accord with what the phenomenon described in (5.1) suggests.

**(5.5) Theorem** (Katz–Klemm–Vafa formula, [36])**.** *The invariants $r^p_{g,m}$ do not depend on the divisibility index, meaning that one has*

$$r^p_{g,m} = r^{m^2p-m^2+1}_{g,1} = R^{m^2p-m^2+1}_g$$

*for all integers $p, m > 0$ and $g \geqslant 0$.*

*They are all determined by the formula*

(5.5.1) $$\sum_{p\geqslant 0}\sum_{g\geqslant 0}(-1)^g\, r^p_g\, \left(y^{\frac{1}{2}}-y^{-\frac{1}{2}}\right)^{2g} q^p = \prod_{n\geqslant 1}\frac{1}{(1-q^n)^{20}(1-yq^n)^2(1-y^{-1}q^n)^2},$$

*where we set $r^p_g := r^p_{g,1}$ (and for convenience, $r^0_0 = 1$, and $r^0_g = 0$ for all $g > 0$).*

Setting $y = 1$ in the formula restricts to the invariants $r^p_0$, and recovers the Yau–Zaslow formula of Theorem (2.1). As a first corollary, one gets that $r^p_g = 0$ if $g > p$, and $r^p_p = (-1)^p(p+1)$. The first values of $r^p_g$ are tabulated below.

| $g$ \ $p$ | 0 | 1 | 2 | 3 | 4 |
|---|---|---|---|---|---|
| 0 | 1 | 24 | 324 | 3200 | 25650 |
| 1 |   | −2 | −54 | −800 | −8550 |
| 2 |   |   | 3 | 88 | 1401 |
| 3 |   |   |   | −4 | −126 |
| 4 |   |   |   |   | 5 |

Table 2: First values of $r^p_g$

### 5.3 – Further Gromov–Witten integrals

We close this section with a short discussion of further results about Gromov–Witten integrals on $K3$ surfaces. They all come from [26].

**(5.6) Hodge integrals with point insertions.** As a direct generalization of the invariants $R^p_{g,m}$ from (5.3.1), one may consider the integrals

$$R^p_{g,k,m} := \int_{[\overline{M}_{g,k}(S,mL)]^{\mathrm{red}}} (-1)^{g-k}\lambda_{g-k} \cup \mathrm{ev}^*_1(\mathrm{pt}) \cup \cdots \cup \mathrm{ev}^*_k(\mathrm{pt}),$$

where $(S, L)$ is a primitive $K3$ surface of genus $p$, $\overline{M}_{g,k}(S, mL)$ is the moduli space of genus $g$ stable maps with $k$ marked points, and $\lambda_i$ is the $i$-th Chern class of the Hodge bundle $\mathbf{E}_{g,k} \to \overline{M}_{g,k}(S, mL)$ as in (5.2).

For *primitive* classes on $K3$ surfaces, the following formula was proved by Maulik, Pandharipande and Thomas [26, Thm. 3]:

(5.6.1) $$\sum_{g=0}^{+\infty}\sum_{p=0}^{+\infty} R^p_{g,k,1}\, u^{2g-2} q^p =$$

$$q\frac{(2\pi)^{12}}{u^2\Delta(q)} \cdot \exp\left(\sum_{g=1}^{+\infty} u^{2g}\frac{B_{2g}}{g(2g)!}E_g(q)\right) \cdot \left(\sum_{m=1}^{+\infty} q^m \sum_{d\mid m}\frac{m}{d}\left(2\sin d\frac{u}{2}\right)^2\right)^k$$

(notation is as in (3.6)).

Note that in the $k = 0$ case, this contains nothing new with respect to the formula given in Theorem (5.5), as the expression (5.6.1) may be deduced from (5.5.1) using known identities, see [26, § 5.4].

**(5.7) Descendent Gromov–Witten invariants.** So far, we have been essentially concerned with the reduced Gromov–Witten invariants (see (3.3))

$$N_g(S, \beta) = \int_{[\overline{M}_{g,g}(S,\beta)]^{\mathrm{red}}} \mathrm{ev}_1^*(\mathrm{pt}) \cup \cdots \cup \mathrm{ev}_g^*(\mathrm{pt})$$

counting curves in the algebraic class $\beta \in H_2(S, \mathbf{Z})$ on the $K3$ surface $S$ and passing through $g$ general fixed points. It is of course possible to pull-back more general cohomology classes $\gamma_i \in H^*(S, \mathbf{Z})$ by the evaluation maps, thus encoding more general incidence conditions than the passing through a given point (although this is not really relevant for surfaces due to the divisor axiom of Gromov–Witten theory). Beware that when doing so one gets integrals that do depend on the class $\beta$ itself, and not only on its self-intersection and divisibility index, as classes in $H^2(S, \mathbf{Z})$ are not monodromy invariant.

A more sensible generalization is to integrate *descendent classes*. Let $\overline{M}_{g,k}(S, \beta)$ be the moduli space of genus $g$ stable maps with $k$ marked points realizing the class $\beta$, and $\mathrm{ev}_1, \ldots, \mathrm{ev}_k$ the corresponding evaluation maps $\overline{M}_{g,k}(S, \beta) \to S$. For all $i = 1, \ldots, k$, define the *i-th cotangent line bundle* $L_i$ to be the line bundle over $\overline{M}_{g,k}(S, \beta)$ the fibre of which over the point $(f : C \to S, p_1, \ldots, p_k)$ is the **C**-line $\Omega^1_{C,p_i}$. The descendent classes on $\overline{M}_{g,k}(S, \beta)$ are those gotten from the Chern classes of these line bundles.

Let $\psi_i := c_1(L_i) \in H^2\big(\overline{M}_{g,k}(S, \beta), \mathbf{Q}\big)$. For all cohomology classes $\gamma_1, \ldots, \gamma_k \in H^*(S, \mathbf{Z})$ and non-negative integers $n_1, \ldots, n_k$, we set

$$(5.7.1) \quad \big\langle \tau_{n_1}(\gamma_1) \cdots \tau_{n_k}(\gamma_k) \big\rangle_g^{S,\beta} := \int_{[\overline{M}_{g,k}(S,\beta)]^{\mathrm{red}}} \psi_1^{k_1} \cup \mathrm{ev}_1^*(\gamma_1) \cup \cdots \cup \psi_k^{k_k} \cup \mathrm{ev}_k^*(\gamma_k)$$

whenever the degree of the integrand equals the (real) dimension $2g + 2k$ of the reduced virtual class, and $\big\langle \tau_{n_1}(\gamma_1) \cdots \tau_{n_k}(\gamma_k) \big\rangle_g^{S,\beta} := 0$ otherwise. These integrals are called *reduced descendent Gromov–Witten invariants*.

The geometric interpretation of the insertion of the classes $\psi_i$ is not straightforward; I refer to [34] and [16] for some discussions about this. See however [16, Thm. 2.2.6], where descendent classes are used to define Gromov–Witten invariants of a projective manifold $X$ *relative* to a smooth very ample hypersurface $Y$, i.e., invariants virtually counting curves in $X$ with prescribed tangency conditions along $Y$.

**(5.8) Quasi-modularity.** The integrals (5.7.1) for fixed integrand and fixed $g$ and divisibility index of $\beta$ are expected to fit together as the Fourier coefficients of a quasi-modular form, as in Theorem (3.5). Due to their dependency on the class $\beta$ and not only on its numerical characters, this is formulated as follows.

Let $S$ be an arbitrarily fixed $K3$ surface possessing an elliptic fibration $\pi : S \to \mathbf{P}^1$ and a section $E$ of $\pi$. Call $\mathsf{e}, \mathsf{f} \in H_2(S, \mathbf{Z})$ the classes of $E$ and the fibres of $\pi$ respectively. It follows from deformation invariance and the same standard degeneration argument as in the proof of Theorem (3.5) that any integral of the form (5.7.1) on an algebraic $K3$ surface equals an integral of the same kind on $S$ with $\beta = a\mathsf{e} + b\mathsf{f}$ for some non-negative integers $a, b$. For all integers $g \geqslant 0$ and $m > 0$, we set

$$F^S_{g,m}\big(\tau_{n_1}(\gamma_1) \cdots \tau_{n_k}(\gamma_k)\big) := \sum_{n \geqslant 0} \big\langle \tau_{n_1}(\gamma_1) \cdots \tau_{n_k}(\gamma_k) \big\rangle_g^{S, m\mathsf{e}+n\mathsf{f}} q^{m(n-m)}$$

as a formal power series in the variable $q$.

**(5.9) Conjecture** (Maulik and Pandharipande, [33, Conj. 3] and [26, § 7.5])**.** *The power series $F^S_{g,m}\big(\tau_{n_1}(\gamma_1) \cdots \tau_{n_k}(\gamma_k)\big)$ is the Fourier expansion in $q$ of a quasi-modular form of level $m^2$ with pole at $q = 0$ of order at most $m^2$.*

(A quasi-modular form of level $N$ with possible pole at $q = 0$ is by definition an element of the **C**-algebra generated by the Eisenstein series $G_1$ (see (3.6)) and modular forms of level $N$; recall in addition that a modular form of level $N$ is a form satisfying the modular equation for transformations in the congruence subgroup $\Gamma_0(N)$ consisting of elements of $\mathrm{PSL}_2(\mathbf{Z})$ congruent to the identity matrix modulo $N$).

For $m = 1$, i.e., for primitive classes, this was proved by Maulik, Pandharipande and Thomas [26, Thm. 4]. Note however that, even in the primitive case, there is no general explicit formula for the modular form in question, as far as I know. Theorem (3.5) provides particular instances of such a formula. At any rate, modularity strongly constrains the invariants and in favorable cases enables one to compute them all (see (6.11) for an example in a different context).

We will only say the following about the proofs of the results presented in this section. A fundamental ingredient for them is the use of other counting invariants than those coming from Gromov–Witten theory, together with correspondence theorems between the two. They are more algebraic in nature than Gromov–Witten invariants, and more agile to study the problems we have been discussing. These invariants virtually count *stable pairs*; they were defined by Pandharipande and Thomas, specifically for threefolds up to now. See [32] for a presentation.

## 6 – Noether–Lefschetz theory and applications

### 6.1 – Lattice polarized $K3$ surfaces and Noether–Lefschetz theory

In this subsection we define Noether–Lefschetz divisors in the moduli spaces of lattice polarized $K3$ surfaces. While the version we will use is the refined one of (6.4), the elementary version of (6.3) is needed to give a proper definition.

Let $\mathbf{L}_{K3} := U^{\oplus 3} \oplus E_8(-1)^{\oplus 2}$ be the $K3$ lattice (see, e.g., [2]). Throughout this subsection, we consider a fixed lattice $\Lambda$ of rank $r$ and signature $(1, r-1)$, together with a primitive embedding $\iota : \Lambda \hookrightarrow \mathbf{L}_{K3}$ (the embedding is *primitive* if the corresponding quotient $\mathbf{L}_{K3}/i(\Lambda)$ is torsion-free).

**(6.1) Definition.** *Let $S$ be a $K3$ surface. A $\Lambda$-polarization on $S$ is a primitive embedding $j : \Lambda \hookrightarrow \mathrm{Pic}(S)$ such that*
  *(i) there is a nef and big class in $j(\Lambda) \subseteq \mathrm{Pic}(S)$;*
  *(ii) there exists an isometry $\phi : H^2(S, \mathbf{Z}) \to \mathbf{L}_{K3}$ such that $\phi \circ j = \iota$.*

*A $\Lambda$-polarized $K3$ surface is a pair $(S, j)$ where $S$ is a $K3$ surface and $j$ is a $\Lambda$-polarization on $S$.*

There exists a moduli space $\mathcal{K}_\Lambda$ of $\Lambda$-polarized $K3$ surfaces, which may be constructed relying on the global Torelli theorem by adapting the method of [37, Exp. XIII, §3].

**(6.2)** Define the discriminant of a rank $s$ lattice $L$ to be the signed determinant

$$\mathrm{Disc}(L) := (-1)^{s-1} \det\bigl(\langle v_i, v_j \rangle\bigr)_{1 \leqslant i, j \leqslant s}$$

where $(v_1, \ldots, v_s)$ is an integral basis of $L$ (the sign has been added to the usual definition so that $\mathrm{Disc}(\Lambda) > 0$); this does not depend on the choice of the basis.

Let $\mathbf{L}$ be a rank $r+1$ lattice with an even symmetric bilinear form, together with a primitive embedding $i : \Lambda \hookrightarrow \mathbf{L}$. There is an invariant of the pair $(\mathbf{L}, i)$ called the coset, which is defined as follows. Consider any vector $v \in \mathbf{L}$ such that $\mathbf{L} = i(\Lambda) \oplus v$; the pairing with $v$ determines an element $\ell_v \in \Lambda^\vee$ in the lattice dual to $\Lambda$. On the other hand, let $G_\Lambda := \Lambda^\vee/\Lambda$ be the quotient of the injection $\Lambda \hookrightarrow \Lambda^\vee$ defined by the pairing on $\Lambda$; it is an abelian group of order $\mathrm{Disc}(\Lambda)$. The *coset* $\delta$ of $(\mathbf{L}, i)$ is the class of $\ell_v$ in $G_\Lambda/\pm$; it does not depend on the choice of $v$.

Two pairs $(\mathbf{L}, i)$ and $(\mathbf{L}', i')$ as above are isomorphic (i.e., there exists an isometry $\phi : \mathbf{L} \to \mathbf{L}'$ such that $\phi \circ i = i'$) if and only if the following two conditions hold: (i) $\mathrm{Disc}(\mathbf{L}) = \mathrm{Disc}(\mathbf{L}')$, and (ii) $\delta(\mathbf{L}, i) = \delta(\mathbf{L}', i')$.

**(6.3) Elementary Noether–Lefschetz divisors.** The *Noether–Lefschetz divisor* $P^\Lambda_{\Delta,\delta} \subseteq \mathcal{K}_\Lambda$ is the closure of the locus of $\Lambda$-polarized $K3$ surfaces $(S, j)$ such that $\mathrm{Pic}(S)$ has rank $r + 1$ and discriminant $\Delta$, and the coset $\delta(\mathrm{Pic}(S), j)$ equals $\delta$.

It follows from the Hodge index theorem that the divisor $P^\Lambda_{\Delta,\delta}$ is empty when $\Delta \leqslant 0$.

**(6.4) Refined Noether–Lefschetz divisors.** Let us now fix an integral basis $\mathbf{v}_\Lambda = (v_1, \ldots, v_r)$ for $\Lambda$, and let $m \in \mathbf{Z}_{>0}$, and $(p, \mathbf{d}) = (p, d_1, \ldots, d_r) \in \mathbf{Z}^{r+1}$. We want to define a Noether–Lefschetz divisor $D^{\mathbf{v}_\Lambda}_{m,p,\mathbf{d}} \subseteq \mathcal{K}_\Lambda$ corresponding to $\Lambda$-polarized $K3$ surfaces $(S, j)$ with an extra class $\beta \in \mathrm{Pic}(S)$ of divisibility index $m$ such that $\langle \beta, \beta \rangle = 2p - 2$, and $\langle \beta, v_i \rangle = d_i$ for all $i = 1, \ldots, r$. The definition depends on the following cases: let

$$\Delta^{\mathbf{v}_\Lambda}_{p,\mathbf{d}} := (-1)^r \begin{vmatrix} \langle v_1, v_1 \rangle & \cdots & \langle v_1, v_r \rangle & d_1 \\ \vdots & \ddots & \vdots & \vdots \\ \langle v_r, v_1 \rangle & \cdots & \langle v_r, v_r \rangle & d_r \\ d_1 & \cdots & d_r & 2p - 2 \end{vmatrix};$$

— if $\Delta^{\mathbf{v}_\Lambda}_{p,\mathbf{d}} > 0$, set

$$D^{\mathbf{v}_\Lambda}_{m,p,\mathbf{d}} := \sum_{\Delta, \delta} \mu^{\mathbf{v}_\Lambda}_{m,p,\mathbf{d}}(\Delta, \delta) \cdot P^\Lambda_{\Delta,\delta}$$

where the sum runs over all $\Delta$ and $\delta$ such that there exists a pair $(\mathbf{L}, i)$ as in (6.2) with $\mathrm{Disc}(\mathbf{L}) = \Delta$ and $\delta(\mathbf{L}, i) = \delta$ (the pair $(\mathbf{L}, i)$ is then unique up to isomorphism), and $\mu^{\mathbf{v}_\Lambda}_{m,p,\mathbf{d}}(\Delta, \delta)$ is the number of elements $\beta \in \mathbf{L}$ having divisibility index $m$ and satisfying $\langle \beta, \beta \rangle = 2p - 2$ and $\langle \beta, v_i \rangle = d_i$ for all $i = 1, \ldots, r$. Note that $\mu^{\mathbf{v}_\Lambda}_{m,p,\mathbf{d}}(\Delta, \delta)$ may be 0; its non-vanishing implies that $\Delta$ divides $\Delta^{\mathbf{v}_\Lambda}_{p,\mathbf{d}}$, so the above sum has only finitely many terms. The condition $\Delta^{\mathbf{v}_\Lambda}_{p,\mathbf{d}} > 0$ implies that any $\beta$ such that $\langle \beta, \beta \rangle = 2p - 2$ and $\langle \beta, v_i \rangle = d_i$ for all $i$ does not belong to $i(\Lambda)$;
— if $\Delta^{\mathbf{v}_\Lambda}_{p,\mathbf{d}} < 0$, set $D^{\mathbf{v}_\Lambda}_{m,p,\mathbf{d}} := 0$;
— if $\Delta^{\mathbf{v}_\Lambda}_{p,\mathbf{d}} = 0$ and $m = \gcd(d_1, \ldots, g_r)$, let $D^{\mathbf{v}_\Lambda}_{m,p,\mathbf{d}}$ be the divisor associated to the dual of the Hodge line bundle $\mathcal{E} \to \mathcal{K}_\Lambda$ (the fibre of $\mathcal{E}$ over the point $(S, i)$ is $H^{2,0}(S)$);
— if $\Delta^{\mathbf{v}_\Lambda}_{p,\mathbf{d}} = 0$ and $m \neq \gcd(d_1, \ldots, g_r)$, set $D^{\mathbf{v}_\Lambda}_{m,p,\mathbf{d}} := 0$.

## 6.2 – Invariants of families of lattice polarized $K3$ surfaces

**(6.5) Families of lattice polarized $K3$ surfaces.** Let $\iota : \Lambda \hookrightarrow \mathbf{L}_{K3}$ be a primitive embedding of a lattice $\Lambda$ of rank $r$ and signature $(1, r-1)$. A *1-parameter family of $\Lambda$-polarized $K3$ surfaces* is a smooth family $\pi : X \to C$ of $K3$ surfaces equipped with line bundles $L_1, \ldots, L_r$ on $X$ such that:
  (i) $X$ is a compact 3-dimensional complex manifold (not necessarily algebraic), $C$ is a complete smooth complex curve, and $\pi$ is a holomorphic submersion;
  (ii) for all $t \in C$, the fibre $X_t$ of $\pi$ over $t$ is a (smooth) $K3$ surface;
  (iii) there exists a linear combination $L^\pi$ of the holomorphic line bundles $L_i$ on $X$, the restriction of which to every fibre of $\pi$ is nef and big;
  (iv) there exists an integral basis $(v_1, \ldots, v_r)$ of $\Lambda$ such that for all $t \in C$, the map $j_t : \Lambda \to \mathrm{Pic}X_t$ defined by $v_i \mapsto L_{i,t}$ (the restriction of $L_i$ to $X_t$) is a $\Lambda$-polarization of $X_t$.

For the remainder of this subsection, we consider a 1-parameter family of $\Lambda$-polarized $K3$ surfaces $(\pi : X \to C, L_1, \ldots, L_r)$ as in (6.5) above.

**(6.6) Noether–Lefschetz numbers.** Let $m \in \mathbf{Z}_{>0}$ and $(p, \mathbf{d}) = (p, d_1, \ldots, d_r) \in \mathbf{Z}^{r+1}$. The *Noether–Lefschetz number* $\mathrm{NL}^\pi_{m,p,\mathbf{d}}$ is defined as

$$\mathrm{NL}^\pi_{m,p,\mathbf{d}} := \int_C f_\pi^* \big( D^{\mathbf{v}_\Lambda}_{m,p,\mathbf{d}} \big),$$

where $f : C \to \mathcal{K}_\Lambda$ is the morphism corresponding to $(\pi : X \to C, L_1, \ldots, L_r)$ by the universal property of $\mathcal{K}_\Lambda$, and $\mathbf{v}_\Lambda$ is the integral basis of $\Lambda$ defined by $(L_1, \ldots, L_r)$ via point (iv) of (6.5).

Note that this is a classical intersection product (i.e., there is no need to define a virtual class), although it may be given by an excess formula in case the image $f_\pi(C)$ is fully contained in the divisor $D^{\mathbf{v}_\Lambda}_{m,p,\mathbf{d}}$.

**(6.7) Gromov–Witten invariants for vertical curve classes.** Although it may not be a projective variety, the total space $X$ carries a $(1,1)$-form $\omega_\pi$ which is Kähler on the fibres of $\pi$; this is sufficient to define a Gromov–Witten theory for non-zero *vertical* classes $\gamma \in H_2(X, \mathbf{Z})^\pi$, i.e., classes $\gamma \in H_2(X, \mathbf{Z})$ such that $\pi_*(\gamma) = 0$ (see [28, §2.1] for details).

In this theory, the moduli spaces of genus $g$ stable maps $\overline{M}_g(X, \gamma)$ all have virtual dimension 0, and we have a set of invariants

$$N^X_{g,\gamma} := \int_{[\overline{M}_g(X,\gamma)]^{\mathrm{vir}}} 1$$

for non-zero vertical classes $\gamma$. We consider the invariants $n^X_{g,\gamma}$ obtained from the $N^X_{g,\gamma}$ by applying the BPS corrections packaged in the formula of (5.4): we let

$$F^X(u,v) := \sum_{g \geqslant 0} \sum_{0 \neq \gamma \in H_2(Z, \mathbf{Z})^\pi} N^X_{g,\gamma} u^{2g-2} v^\gamma$$

as a formal power series in the variables $u, v$, where the powers of $v$ are indexed by $H_2(Z, \mathbf{Z})^\pi$, and set

$$F^X(u,v) := \sum_{g \geqslant 0} \sum_{0 \neq \gamma \in H_2(Z, \mathbf{Z})^\pi} n^X_{g,\gamma} u^{2g-2} \left( \sum_{d > 0} \frac{1}{d} \left( \frac{\sin d\frac{u}{2}}{\frac{u}{2}} \right)^{2g-2} v^{d\gamma} \right).$$

Finally, for a non-zero multidegree $\mathbf{d} = (d_1, \ldots, d_r) \in \mathbf{Z}^r$, we let $n^X_{g,\mathbf{d}}$ be the invariant counting genus $g$ stable maps in vertical classes of degree $d_1, \ldots, d_r$ with respect to $L_1, \ldots, L_r$ respectively, i.e.,

(6.7.1) $$n^X_{g,\mathbf{d}} := \sum_{\gamma \in H_2(X, \mathbf{Z})^\pi : \int_\gamma L_i = d_i} n^X_{g,\gamma}.$$

**(6.8) Reduced Gromov–Witten invariants of $K3$ fibres.** We also consider the invariants $r^p_{g,m}$ for $K3$ surfaces which have been defined in (5.4); recall that they are the reduced Hodge integrals (5.3.1) put in BPS form. Note that we need to maintain the dependency on the divisibility index $m$, because Theorem (6.9) below is needed for the proof of the independence on $m$ conjectured by Yau–Zaslow.

A multidegree $\mathbf{d} = (d_1, \ldots, d_r) \in \mathbf{Z}^r$ is positive with respect to $L^\pi$ if for any line bundle $M$ on some fibre $X_t$ of $\pi$, $(M, L_{i,t}) = d_i$ for all $i$ implies $(M, L^\pi) > 0$; since $L^\pi$ is a linear combination of the $L_i$, this is an elementary linear algebraic condition.

**(6.9) Theorem** ([28])**.** *Let* $\mathbf{d} = (d_1, \ldots, d_r) \in \mathbf{Z}^r$ *be a multidegree which is positive with respect to* $L^\pi$. *Then*

$$(6.9.1) \qquad n_{g,\mathbf{d}}^X = \sum_{p=0}^{+\infty} \sum_{m=1}^{+\infty} r_{g,m}^p \cdot \mathrm{NL}_{m,p,\mathbf{d}}^\pi.$$

(This is stated in [28] in the $r = 1$ case (*i.e.*, $\mathbf{d} \in \mathbf{Z}$), but as noted in [23] the same proof goes through in general).

The philosophy behind this relation is rather natural, and ought to be compared to the discussion of Section 5.1 above. Consider the genus $g = 0$ case for simplicity; then the invariant $n_{0,\mathbf{d}}^X$ counts vertical rational curves in $X$ of prescribed degrees with respect to $L_1, \ldots, L_r$, and these are virtually in finite number. There are on the other hand finitely many members of the family $\pi$ with algebraic divisor classes of the prescribed degrees with respect to $L_1, \ldots, L_r$, each of which provides a finite number of rational curves. The theorem morally says that the number of rational curves in $X$ is the sum of these isolated contributions from the fibres.

Of course the actual story is more complicated than this, if only because of the existence of 1-dimensional families of rational curves on $X$, coming from finitely many rational curves in *all* $K3$ members of the family (which all have algebraic divisor classes, as they are $\Lambda$-polarized), in spite of the virtual dimension being 0. In other words, a Calabi-Yau threefold $X$ as in Theorem (6.9) above is far from satisfying the same properties than the perturbations of the twistor families of algebraic $K3$ surfaces on which BPS numbers are supposed to count curves.

## 6.3 – Application to the Yau–Zaslow conjecture

In this subsection we give an outline of the proof by Klemm, Maulik, Pandharipande, and Scheidegger of the Yau–Zaslow conjecture (Theorem (4.14) above). Recall that the invariants $r_{g,m}^p$ being invariant under algebraic deformations of the $K3$ surface, it is enough to prove the result for our favourite $K3$ surface. These invariants for certain elliptic $K3$ surfaces are expressed by means of the relation (6.9.1) for a particular family.

**(6.10) The *STU* model.** The central character of the proof is a smooth projective Calabi–Yau 3-fold $X$, known as the *STU* model and coming from physics (quoting [23]: the letter $S$ stands for the dilaton and $T$ and $U$ label the torus moduli in the heterotic string). It is contructed as an anticanonical section of a smooth projective toric 4-fold $Y$ defined by an explicit fan in $\mathbf{Z}^4$.

The variety $X$ has the structure of a fibration $\pi : X \to \mathbf{P}^1$, the general fibre of which is a smooth $K3$ surface, itself with an elliptic fibration. It comes with two line bundles $L_1, L_2 \to X$, defining a $\Lambda$-polarization on the family $\pi : X \to \mathbf{P}^1$ (leaving aside the fact that there are inevitably singular members), where $\Lambda$ is the lattice with intersection form

$$\begin{pmatrix} 0 & 1 \\ 1 & 0 \end{pmatrix}.$$

The family $\pi$ has the shape of a Lefschetz pencil, in particular each of its singular members has a unique ordinary double point as its only singularity. One may thus build an actual family $\tilde{\pi} : \tilde{X} \to C$ of $\Lambda$-polarized $K3$ surfaces from $\pi$ as follows. One first performs a base change by $t \mapsto t^2$ around each singular member; to do so, one considers the $2:1$ covering $\varepsilon : C \to \mathbf{P}^1$ with branch divisor $\mathrm{Disc}(\pi)$, the set of points above which $\pi$ fails to be smooth, and let $\pi^\flat : X^\flat \to C$ be the family obtained from $\pi$ by applying the base change $\varepsilon : C \to \mathbf{P}^1$. The new total space $X^\flat$ is singular, precisely it has an ordinary 3-fold double point at each singular point of a fibre (analytically locally around such a point, $X$ is defined by the equation $x^2 + y^2 + z^2 = t$ in a 4-dimensional complex ball, hence $X^\flat$ is defined by $x^2 + y^2 + z^2 = t^2$). One then chooses for $\tilde{X}$ any small resolution of all these singularities: this may be understood as first blowing-up once all singular points, and then contracting one ruling of each exceptional divisor (they are all smooth quadric surfaces). This has the effect of replacing each fibre of $X^\flat$ by its minimal model.

One may determine the number of singular members of $\pi$ by the same topological Euler characteristic computation as in Section 2.1. The Euler number $e(X)$ is found to be $-480$ by toric intersection computations in the 4-fold $Y$, hence the number of singular fibres equals

$$e(K3) \cdot e(\mathbf{P}^1) - e(X) = 528$$

($e(K3)$ denotes the common Euler number of all $K3$ surfaces, i.e., 24).

**(6.11) Modularity for Noether–Lefschetz numbers.** It is a stunning application of a theory developed by Borcherds and Kudla–Millson (see [28, 23] and the references therein) that the Noether–Lefschetz numbers of a family of $\Lambda$-polarized $K3$ surfaces fit into a vector valued modular form.

In order to state this precisely, let us return to the setup of (6.6) for a moment (see [28, § 4] for a complete treatment). One may define divisors $D_{p,\mathbf{d}}^{\mathbf{v}_\Lambda}$ and numbers $\mathrm{NL}_{p,\mathbf{d}}^\pi$ by dropping the requirement on the divisibility index $m$ in (6.4). It is an elementary result [23, Lemma 1] that the full set of the numbers $\mathrm{NL}_{p,\mathbf{d}}^\pi$ determine the refined Noether–Lefschetz numbers $\mathrm{NL}_{m,p,\mathbf{d}}^\pi$. Let $\mathrm{Mp}_2(\mathbf{Z})$ be the metaplectic double cover of $\mathrm{SL}_2(\mathbf{Z})$. There is a canonical representation

$$\rho_\Lambda^* : \mathrm{Mp}_2(\mathbf{Z}) \to \mathrm{End}\bigl(\mathbf{C}[G_\Lambda]\bigr)$$

associated to $\Lambda$ (recall that $G_\Lambda = \Lambda^\vee/\Lambda$).

**(6.12) Theorem** (Borcherds, Kudla–Millson, Maulik–Pandharipande)**.** *There exists a vector-valued modular form*

$$\Phi^\pi(q) = \sum_{\gamma \in G_\Lambda} \Phi_\gamma^\pi(q)\, u^\gamma \in \mathbf{C}[[q^{\frac{1}{2\mathrm{Disc}\Lambda}}]] \otimes \mathbf{C}[G_\Lambda]$$

*of weight $\frac{22-r}{2}$ and type $\rho_\Lambda^*$, such that the Noether–Lefschetz number $\mathrm{NL}_{p,\mathbf{d}}^\pi$ is the coefficient of $\Phi_\gamma^\pi$ in $q$ to the power $\frac{\Delta_{p,\mathbf{d}}^{\mathbf{v}_\Lambda}}{2\mathrm{Disc}\Lambda}$, where $\gamma \in G_\Lambda$ is any of the two liftings of the coset $\delta_{p,\mathbf{d}}^{\mathbf{v}_\Lambda} \in G_\Lambda/\pm$ represented by the linear functional $v_i \mapsto d_i$.*

Taking advantage of the strong structure results for modular forms, Maulik and Pandharipande are able to use this theorem to derive explicitly the Noether–Lefchetz numbers of classical families of $K3$ surfaces of genus $2 \leqslant p \leqslant 5$ (i.e., double planes and complete intersection $K3$'s). A similar calculation is carried out in [23] for the *STU* family, as one of the key steps in the proof of the Yau–Zaslow conjecture.

**(6.13)** Let us now return to the setup of (6.10). Theorem (6.12) tells that the Noether–Lefschetz numbers of the family $\tilde\pi : \tilde{X} \to C$ are the Fourier coefficients of a scalar modular form of weight 10. The vector space of such forms has dimension 1 and is generated by the Eisenstein series

$$E_5(q) = E_2(q) E_3(q) = 1 - 264 \sum_{n=1}^{+\infty} \sigma_9(n) q^n$$

[39, § VII.3.2] (notation as in (3.6)). It follows that it is enough to know one Noether–Lefschetz number to determine the full modular form, and since we do know one of them, namely the number 528 of singular members of the *STU* family, one obtains that the number $\mathrm{NL}^{\tilde{\pi}}_{p,d_1,d_2}$ is the coefficient in $q$ to the power $\frac{1}{2}\Delta(p,d_1,d_2)$ of the modular form $-4E_2(q)E_3(q)$, where

$$\Delta(p,d_1,d_2) = \begin{vmatrix} 0 & 1 & d_1 \\ 1 & 0 & d_2 \\ d_1 & d_2 & 2p-2 \end{vmatrix}.$$

**(6.14) Mirror symmetry.** The *STU* model $X$ being an anticanonical section of a smooth semi-positive toric variety, its genus 0 Gromov–Witten invariants are known by mathematically proven mirror symmetry results. This gives the corresponding invariants of $\tilde{X}$, the latter being twice those of $X$ [28, § 5.2].

Precisely, Givental has proven the relation of the genus 0 Gromov–Witten invariants of $X$ by mirror transformation to hypergeometric solutions of the Picard–Fuchs equations of the Batyrev–Borisov mirror, see [28, 23] and the references therein. This gives the following formula of Klemm–Mayr–Lerche [23, Prop. 5]

$$(6.14.1) \qquad \sum_{(d_1,d_2)\in\mathcal{P}} (d_2)^3 N^X_{0,(d_1,d_2)} q_1^{d_1} q_2^{d_2} = -2 + 2\frac{E_2(q_1)E_3(q_1)}{(2\pi)^{-12}\Delta(q_1)} \frac{E_2(q_2)}{j(q_1) - j(q_2)},$$

where

$$j(q) := 1728 \frac{(60G_2(q))^3}{\Delta(q)} = (2\pi)^{12} \frac{E_2(q)^3}{\Delta(q)} = \frac{1}{q} + 744 + 196884q + \cdots$$

(notation as in (3.6)) is the normalized $j$ function, $\mathcal{P} = \{(d_1,d_2) \neq (0,0) : d_1 \geqslant 0, d_1 \geqslant -d_2\}$, and $N^X_{0,(d_1,d_2)}$ is defined by formula (6.7.1) from the various $N^X_{0,\gamma}$, $\gamma \in H_2(X,\mathbf{Z})^\pi$.

**(6.15) Conclusion: the Harvey–Moore identity.** Using the fact that the lattice $\Lambda$ has rank 2, Klemm–Maulik–Pandharipande–Scheidegger then show that the invariants $r^p_{0,m}$ are uniquely determined by the relations (6.9.1) for the family $\tilde{\pi} : \tilde{X} \to C$, and the numbers $n^{\tilde{X}}_{0,(d_1,d_2)}$ and $\mathrm{NL}^{\tilde{\pi}}_{m,p,(d_1,d_2)}$ [23, Prop. 3]. The two latter sets of numbers being known by the results of (6.11) and (6.14), it is therefore enough, in order to end the proof of Theorem (4.14), to show that the numbers $r^p_{0,m}$ predicted by the Yau–Zaslow conjecture (i.e., $r^p_{0,m} = r^{m^2p-m^2+1}_{0,1}$ together with the formula of (2.1) giving the $r^p_{0,1}$'s) indeed fit in the relations (6.9.1). This takes the form of an identity between modular forms: let

$$f(z) := \frac{E_2(z)E_3(z)}{(2\pi)^{-12}\Delta(z)} = \sum_{n=-1}^{+\infty} c(n)q^n, \quad q = e^{2\pi i z};$$

what has to be proven is

$$\frac{f(z_1)E_2(z_2)}{j(z_1) - j(z_2)} = \frac{q_1}{q_1 - q_2} + E_2(z_2) - \sum_{d,k,l>0} l^3 c(kl) q_1^{kd} q_2^{ld}.$$

This is the Harvey–Moore identity, which has been proven by Zagier, see [23, § 4.2].

## References

**Chapters in this Volume**[7]

[**] Th. Dedieu and E. Sernesi, *Deformations of curves on surfaces.*

---